\documentclass[11pt]{article}
\usepackage{amsfonts}
\usepackage{amsmath}
\usepackage{graphicx}
\usepackage{bbm}
\usepackage[textsize=small]{todonotes}
\newtheorem{theorem}{Theorem}[section]

\newtheorem{assumption}[theorem]{Assumption}

\newtheorem{lemma}[theorem]{Lemma}

\newtheorem{proposition}[theorem]{Proposition}

\newtheorem{preremark}{Remark}[section]
  \newenvironment{remark}%
    {\begin{preremark}\rm}{\end{preremark}}

\newcommand{\noi}{\noindent}
\newcommand{\RR}{\mathbb{R}}
\newcommand{\NN}{\mathbb{N}}
\newcommand{\EE}{\mathbb{E}}
\newcommand{\PP}{\mathbb{P}}

\newcommand{\clb}{\mathcal{B}}
\newcommand{\clu}{\mathcal{U}}
\newcommand{\clf}{\mathcal{F}}

\newcommand{\cls}{\mathcal{S}}
\newcommand{\cle}{\mathcal{E}}
\newcommand{\bfr}{\boldsymbol{R}}

\begin{document}
\title{\textbf{Admission Control for Multidimensional Workload  with Heavy Tails and Fractional Ornstein-Uhlenbeck Process}}
\author{Amarjit Budhiraja\thanks{This research is partially supported by the National Science Foundation
	(DMS-1004418, DMS-1016441), the Army Research Office (W911NF-10-1-0158) and the US-Israel Binational Science Foundation (Grant 2008466).}, Vladas Pipiras and Xiaoming Song}

\date{\today}
\maketitle
\begin{abstract}
The infinite source Poisson arrival model with heavy-tailed workload distributions has attracted much attention, especially in the modeling of data packet traffic in communication networks. 
In particular, it is well known that under suitable assumptions on the source arrival rate, the centered and scaled cumulative workload process for the underlying processing system can be approximated by fractional Brownian motion. 
In many applications one is interested in the stabilization of the work inflow to the system by modifying the net input rate, using an appropriate admission control policy.  In this work we study a  natural family of admission control policies which keep the associated scaled cumulative workload asymptotically close to a pre-specified linear trajectory, uniformly over time.
 Under such admission control policies and with natural assumptions on  arrival distributions, suitably scaled and  centered cumulative workload processes are shown to converge weakly in the path space to the solution of a $d$-dimensional stochastic differential equation (SDE) 
driven by a  Gaussian process. It is shown that the admission control policy achieves moment stabilization in that the second moment of the solution to the SDE (averaged over the $d$-stations) is bounded uniformly for all times. In one special case of  control policies, as time approaches infinity, we obtain a 
fractional version of  a stationary Ornstein-Uhlenbeck process that is driven by fractional Brownian motion with Hurst parameter $H>\frac{1}{2}$.\\

\noi {\bf AMS 2010 subject classifications:}
60G15, 60G18, 60G57, 90B15.

\noi {\bf Keywords:}
Poisson random measures, Gaussian random measures, heavy-tailed distribution, self-similarity, fractional Brownian motion, fractional Ornstein-Uhlenbeck process, admission control.

\end{abstract} 
\section{Introduction } \label{sec:sec1}

This work is motivated by a general workload model for data traffic in communication networks considered by Kurtz \cite{K1}. In this model, a large number of sources input work into a system. Let $N(t)$ be the number of source activations up to time $t$. For the activation of the $i$th source, let $X_i(s)$ denote the cumulative work input into the system during the first $s$ units of time that the source is on, and let $\tau_i$ denote the length of time that the $i$th source remains active ($i$th session length). Then, the total work input into the system up to time $t$ is given by
\begin{equation}\label{U}
W(t)=\int_0^tX_{N(s)}(\tau_{N(s)}\wedge(t-s))dN(s).
\end{equation} 
It is assumed that $(X_i, \tau_i)$ are i.i.d.\ with distribution $\nu$ on $D_{\RR_+}[0,\infty)\times[0,\infty)$, where $D_{\RR_+}[0,\infty)$ denotes the space of $\RR_+$-valued right-continuous functions on $[0,\infty)$ having left limits (RCLL).

The paper \cite{K1} assumes that the arrival (counting) process $N(t)$ is characterized by an intensity $\lambda$, which may depend on the past of $N$, $W$ and another process $Q$, representing the number of active sources. For simplicity and since only this case will be considered below, assume that the intensity $\lambda$ depends only on the total workload $W(t)$ and time $t$ as
\begin{equation}
\lambda=\lambda(t,W(t)), 
\end{equation}
and focus only on the process $(N,W)$. Letting $\xi$ be a Poisson random measure on $[0,\infty)\times D_{\RR}[0,\infty)\times[0,\infty)$ with intensity measure $\eta=m\times\nu$ ($m$ being Lebesgue measure), the process $(N,W)$ can be represented as
\begin{align}\label{N-U-1}
N(t)=&\  \xi(B(t)),\nonumber\\
W(t)=&\int_{B(t)}u(r\wedge (t-\gamma(s)))\xi(ds,du,dr),
\end{align}
where
\begin{equation}\label{N-U-2}
B(t)=\{(s,u,r): s\leq \int_0^t\lambda(z,W(z))dz\}
\end{equation}
and $\gamma(s)$ satisfies
\begin{equation}\label{def-gamma}
\int_0^{\gamma(s)}\lambda(z,W(z))dz=s.
\end{equation}
See Section \ref{sec:sec2.1} for some discussion on this system of equations. 
The key results of  \cite{K1} are the law of large numbers and the central limit theorem for the scaled system $(X_n,Y_n)$, where
\begin{equation}\label{N-U-n}
X_n(t)=\frac{1}{n}N_n(t),\, Y_n(t)=\frac{1}{n}W_n(t),
\end{equation}
and $(N_n, W_n)$ are defined as in \eqref{N-U-1} but using a Poisson random measure $\xi_n$ with intensity measure $nm\times\nu$ and an intensity $\lambda_n(t,w(t))=n\lambda(t, n^{-1}w(t))$. Under suitable assumptions, $(X_n,Y_n)$ converges in probability to $(X,Y)$ satisfying 
\begin{equation}
X(t)=\int_0^t\lambda(s,Y(s))ds,\, Y(t)=\int_0^t\mu(t-r)\lambda(s,Y(s))ds,
\end{equation}
where $\mu(t)=\EE(X_i(\tau_i\wedge t))$ (Theorem 2.1 in \cite{K1}). Under suitable assumptions, the scaled and centered process $(\tilde{X}_n,\tilde{Y}_n)=\sqrt{n}(X_n-X,Y_n-Y)$ converges in distribution to $(\tilde{X},\tilde{Y})$ satisfying
\begin{align}\label{limit}
\tilde{X}(t)=&\ \Xi(B(t))+\int_0^t\lambda_y(s,Y(s))\tilde{Y}(s)ds,\nonumber\\
\tilde{Y}(t)=&\int_{B(t)}u(r\wedge(t-\gamma(s)))\Xi(ds,du,dr)+\int_0^t\mu(t-s)\lambda_y(s,Y(s))\tilde{Y}(s)ds,
\end{align} 
where $B(t)$ and $\gamma(s)$ are defined as in \eqref{N-U-2}-\eqref{def-gamma} but replacing $W$ by $Y$, and $\Xi$ is a Gaussian random measure with the control measure
\begin{equation}\label{control-mea}
\EE(\vert \Xi(ds,du,dr)\vert^2)=ds\nu(du,dr)
\end{equation}
(Theorem 2.2 in \cite{K1}).

One special case of $(X_i,\tau_i)$ is particularly interesting in the context of modeling data traffic in modern communication networks. This is the case where
\begin{equation}\label{rate-1}
X_i(s)=s
\end{equation}
and $\tau_i$ are heavy tailed with the distribution
\begin{equation}\label{dist-1}
\nu(dr)=(\beta-1)\theta(\theta r+1)^{-\beta}dr,\, r\geq 0,
\end{equation}
where $\beta\in(2,3)$ is the tail index and $\theta\in(0,\infty)$ is  a scale parameter. (Note that we use the same notation $\nu$ to denote the distribution of $\tau_i$ alone.) When intensity $\lambda$ is constant and \eqref{rate-1}-\eqref{dist-1} are assumed, it is well known that, under suitable assumptions and proper scaling, the cumulative workload process converges to fractional Brownian motion. A version of this fact appears in Section 4 of \cite{K1}. But see also \cite{HRS, KT, KL, LTWW, MRRS,PTL}.

For the convergence to fractional Brownian motion, in a scaled system, it is also necessary to rescale the measure $\nu(dr)$ in \eqref{dist-1}. One way to see this is to observe that without rescaling, the Gaussian random measure $\Xi$ in \eqref{limit}-\eqref{control-mea} is not self-similar in the variable $r$. For this reason (see for example, \cite{K1}, Section 4 for the case when $\lambda$ is constant), it is natural to scale the measures as:
\begin{equation}\label{dist-2}
\nu_n(dr)=(\beta-1)n\theta(n\theta r+1)^{-\beta}dr.
\end{equation}

This case is not included in Theorems 2.1 and 2.2 of \cite{K1} which, although allow for state dependent $\lambda$, treat the scaled system \eqref{N-U-n}  that has no scaling in the intensity measure $\nu$, and hence the key convergence results of \cite{K1} for non-constant $\lambda$ cannot be applied with \eqref{dist-2}. In fact, as already suggested by the result in Section 4 of \cite{K1} for the constant $\lambda$ case, dealing with \eqref{dist-2} for non-constant $\lambda$ is expected to be more involved. For example, a natural normalization in this case is no longer  $\sqrt{n}$.

Models where the arrival intensity $\lambda$ is a function of the state process are natural when one considers control mechanisms for regulating the amount of work in the system.  A common form of a control policy
that aims to appropriately balance long processing delays with low processor utilization, consists of suitably decreasing the input rate when the workload in the system is very high and increasing the rate when it drops too low.  Study of
asymptotic behavior of the workload process with heavy-tailed session length distributions, under such state feedback control mechanisms is the subject of current work.
We shall consider a scaled multidimensional system where the session lengths are distributed according to $\nu_n$ as in  \eqref{dist-2}, and establish limit theorems for settings where   $\lambda$ is state dependent. 
We are particularly interested in the design of control policies that keep the net workload (asymptotically) close to a pre-specified linear trajectory such that the variability (suitably scaled) is bounded uniformly in time.
The slope of the linear trajectory represents the system processing rate and thus such control policies yield uniform in time reliability bounds on probabilities of processor underutilization and overload.

Let us now describe briefly our model and results established below.
We suppose that a system consists of $d$ processing stations, and that workload arrives to each station (independently of others) as before.
The intensity $\lambda$ controlling the arrival rate, however, now depends on the average total workload across all the stations. More specifically, denoting the total cumulative workload at the $i$th station by $y_i(t)$ and their average $\bar{y}(t)$, we suppose that
\begin{equation}
\lambda=f(t,\bar{y}(t)).
\end{equation}
More specifically yet, we will work with special intensities $\lambda$ having the form
\begin{equation}\label{inten}
\lambda=f(t,\bar{y}(t))=\exp\{-g(\bar{y}(t)-bt)\},
\end{equation}
for some $b>0$ and function $g$. The constant $b$ represents the processing rate at each station, although processing of work is not explicitly included in our model
and plays no role in the analysis.
The function $g$ will satisfy the assumption stated next.
\begin{assumption}
	\label{ass:ass1}
	$g(0)=0$. The function $g$ is twice differentiable and its first and second derivatives  $g'$ and $g''$ satisfy 
	\[
	0<\ell\leq g'(x)\leq L,\ \forall x\in\mathbb{R},
	\]
	 and 
	 \[
	 \vert g''(x)\vert\leq L,\ \forall x\in\mathbb{R}
	 \]
	   for some  $\ell, L \in (0, \infty)$.
 \end{assumption}
The above assumption will be taken to hold throughout this work and will not be explicitly noted in the statements of various results.
Note that under this assumption $g$ is a strictly increasing function and $g(u)>0$ if $u>0$ and $g(u)<0$ if $u<0$. From the properties of $g$, we see that the intensity $\lambda$ in \eqref{inten} has a natural physical interpretation: the intensity of session arrivals at the $i$th station increases when $\bar{y}(t)$ drops below $bt$ while it decreases when $ \bar{y}(t)$ exceeds $bt$. We will refer to $g$ as an  \textit{admission control policy}. 

Our scaled system will be characterized by independent Poisson random measures $\xi_{n,i}$ having common intensity measure $n^\alpha m\times\nu_n$ where  $\nu_n$ is as in \eqref{dist-2}
and the cumulative workload process $Y_{n,i}(t)$, unlike \eqref{N-U-n} will now be normalized by a factor of $n^{\alpha-1}$, rather than $n$ (see \eqref{e-2-2}). We will assume that
\begin{equation}
	\label{eq:ab1210}
	\alpha \in \left(\beta -1, \min\{3\beta-5, 5-\beta\}\right).
	\end{equation}
The reason for such choice of $\alpha$, and for the normalization $n^{\alpha-1}$ will be given below (see Remark \ref{Remark-2-2}). 

Precise evolution equations for $Y_{n,i}$ are given in Section \ref{sec:sec2}. 
We now give a brief description of our main results.  In Theorem \ref{conv-prob} we prove a law of large numbers result stating that, as $n \to \infty$,  $\mathbf{Y}_{n}=(Y_{n,1},\dots,Y_{n,d})^T$ converges in
probability, in $D_{\mathbb{R}_+^d}[0,\infty)$
to a continuous (non-random) trajectory $\mathbf{U}=(U,\dots,U)^T$, where $U$ is characterized as the unique solution of an  ordinary differential equation (ODE)  (see \eqref{e-2-4}), and a rate of convergence is given as well.  The solution $U$ has the property that
$\sup_{t\ge 0}|U(t) - bt| < \infty$.  In fact, with a particular choice of $b$, namely $b = \frac{1}{\theta(\beta-2)}$, we have $U(t) = bt$ for all $t$.

Next, we study the fluctuations of $\mathbf{Y}_{n}$.  In Theorem \ref{CLT} we show that suitably centered and normalized form of $\mathbf{Y}_n$, denoted as $\mathbf{Z}_n$ (see \eqref{eq:ab348}), converges in distribution, in 
$D_{\mathbb{R}^d}[0,\infty)$ to the solution $\mathbf{Z} = (Z_1, \dots, Z_d)^T$ of a $d$-dimensional stochastic differential equation (see \eqref{linear}), driven by $d$ independent Gaussian processes $R_i$, $i = 1, \dots , d$.
 The moment stabilization property of the admission control policy is demonstrated in Theorem \ref{momstab}, which says that $\sup_{t\ge 0} \EE|\bar{Z}(t)|^2 < \infty$, where $\bar{Z}=\frac{1}{d}\sum_{i=1}^dZ_i$.

We remark that, in the case when $b = \frac{1}{\theta(\beta-2)}$, one can achieve the law of large number limit of $bt$ by simply taking the admission control policy to be $g\equiv 0$ (this function obviously does not satisfy Assumption \ref{ass:ass1}).  However, in the case when $g\equiv0$, the limit process obtained from the fluctuation central limit theorem will have variance that increases to $\infty$ as $t \to \infty$.  

Finally, we show that  in one particular case, the average of the limit process $\bar{Z}$ is driven by a Gaussian $H$-self-similar process $\bar{R}$ with $H=\frac{4-\beta}{2}>\frac{1}{2}$. The driving process $\bar{R}$ is not  fractional Brownian motion since it does not have stationary increments. This is directly related to the fact that the limit process $\bar{Z}$ satisfies $\bar{Z}(0)=0$ and hence  is not stationary. 
 The process $\bar{Z}(T+ \cdot)$ is expected to become stationary  as $T\to\infty$. Similarly, the driving process $\bar{R}$ is expected to have stationary increments in the long run (i.e. $\bar{R}(T+\cdot)-\bar{R}(T)$ approaches a process with
stationary increments, as $T\to\infty$).  
We study the asymptotic behavior of the process
$\bar{Z}(T+\cdot)$ as $T \to \infty$ in Theorem \ref{ouplim}.
  For simplicity we restrict here to the case $b = \frac{1}{\theta(\beta-2)}$.  It is shown that, as $T \to \infty$, the process
$\bar{Z}(T+\cdot)$ converges in distribution, in $C_{\mathbb{R}}[0,\infty)$ to a stationary Ornstein-Uhlenbeck process driven by fractional Brownian motion with Hurst parameter $H=\frac{4-\beta}{2}>\frac{1}{2}$.

The paper is organized as follows. We state all the results in Section \ref{sec:sec2}. Section \ref{sec:sec3} contains the proofs of Proposition \ref{L-3-1} and Theorem \ref{conv-prob}. The proof of the central limit theorem will be provided in Section \ref{sec:sec4}. In Section \ref{sec:sec5} we represent the limit (centered) station average process $\bar Z$  as an integral with respect to a Gaussian process and give the proof of Theorem \ref{momstab} on the moment stabilization property of the admission control policy $g$. Section 6 is devoted to the study of the
 asymptotic  behavior, as $T \to \infty$, of the  process $\bar Z$ obtained from the central limit theorem, and the proof of Theorem \ref{ouplim} is given.

The following notation will be used.
We denote the set of non-negative integers by $\NN$ and non-negative reals by $\RR_+$.
For a Polish space $S$, $C_S[0, \infty)$ (resp. $D_S[0, \infty)$) will denote the space of continuous (resp. RCLL) functions endowed with the local uniform (resp. Skorohod) topology.
$C$ will denote generic constants in $(0, \infty)$ whose value may change from one proof to next.

\bigskip

\setcounter{equation}{0}
\section{Model Formulation and Main Results}
\label{sec:sec2}
We begin in this section with the evolution equations for the unscaled system.

\subsection{Unscaled System}\label{sec:sec2.1}
  Let $\xi_{0,i}$, $i = 1, \dots, d$, be independent Poisson random measures on $[0,\infty)\times[0,\infty)$ with common intensity $\eta=m\times\nu$, where $m$ denotes the Lebesgue measure on $[0,\infty)$ and $\nu$ is given in \eqref{dist-1}. Then, $\xi_{0,i}$ can be represented as
\[
\xi_{0,i}=\sum_{j=1}^\infty\delta_{(S_{i,j},\tau_{i,j})},
\]
where $0<S_{i,1}<S_{i,2}<\cdots$ are the jump times of independent unit rate  Poisson processes for $i=1,\dots,d$ and $\tau_{i,j}$ are i.i.d.\ with distribution $\nu$. These Poisson random measures will be the building blocks for our counting processes
$N_i$ with desired intensities.

Let $f:\mathbb{R}_+\times\mathbb{R}\to\mathbb{R}_+$ be a function of the form 
\begin{equation}
f(t,y)=\exp\{-g(y-bt)\},
\end{equation}
 where $g:\mathbb{R}\to\mathbb{R}$ is a function satisfying Assumption \ref{ass:ass1}. Let $\mathbf{X}_0 = (X_{0,1}, \dots, X_{0,d})^T$ and $\mathbf{Y}_0 = (Y_{0,1}, \dots, Y_{0,d})^T$ be $\NN^d$ and $\RR_+^d$-valued RCLL processes given through the following system of equations:
%
%
%
\begin{align}
X_{0,i}(t)=&\ N_{0,i}(t) =\ \xi_{0,i}(B_0(t))=\sum_{j=1}^\infty1_{\{S_{i,j}\leq \Lambda_0(t)\}},\nonumber\\
Y_{0,i}(t) =&\int_{B_0(t)}r\wedge (t-\gamma_0(s))\xi_{0,i}(ds,dr) \nonumber\\
=&\sum_{j:S_{i,j}\leq\Lambda_0(t)}\tau_{i,j}\wedge(t-\gamma_0(S_{i,j})),\label{N-U}
\end{align}
where
\begin{align}\label{e-1-2}
\Lambda_0(t) = &\, \int_0^t f(s,\bar{Y}_0(s))ds,\ \bar{Y}_0(t) = \frac{1}{d} \sum_{i=1}^d Y_{0,i}(t),\
 \gamma_0(t) =\ \Lambda_0^{-1}(t),
 \end{align}
 and
\[
B_0(t)= [0, \Lambda_0(t)] \times [0, \infty).
\]
Note that from Assumption \ref{ass:ass1} $\Lambda_0$ is continuous and strictly increasing.  Therefore, $\gamma_0$ is well defined and continuous as well.
$\gamma_0(S_{i,j})$ is the $j$th activation time at the $i$th station, that is, the $j$th jump time of $N_{0,i}(t)$. For $\gamma_0(S_{i,j})\leq t$, $t-\gamma_0(S_{i,j})$ is the 
amount of time up to $t$ since the $j$th session activation at the $i$th station and $T_{i,j}=\gamma_0(S_{i,j})+\tau_{i,j}$ is the end time of the $j$th session at the $i$th station.  Thus $\tau_{i,j}\wedge (t-\gamma_0(S_{i,j}))$ is the work input by the $j$th activated source at
the $i$th station, up to time $t$. 

From Assumption \ref{ass:ass1} it follows that
\begin{equation}\label{new-eq-2-4}
f(t,y) = \exp\{-g(y-bt)\} \le \exp\{-g(-bt)\},\ \forall\ y \ge 0.
\end{equation}
In particular, $f$ is a strictly positive function that is locally bounded, namely $$\sup_{t \in [0,T], y \in \RR_+} f(t,y) < \infty,\; \ \forall\  T>0.$$  From this it follows that there is a unique solution to the system
of equations \eqref{N-U}-\eqref{e-1-2}.  Indeed, the solution can be constructed recursively between successive jump times of the underlying Poisson random measures $\{\xi_{0,i}, i = 1, \dots, d\}$.  To see the basic construction, we consider for simplicity
the case $d=1$.  The general multi-dimensional case can be treated similarly by arranging the jumps $\{S_{i,j}\}$ in an increasing order.  Simplifying notation, denote by $0<S_1<S_2<\dots$ the jump times
of $\xi = \xi_{0,1}$ and denote the corresponding session lengths by $\tau_1, \tau_2, \ldots$.  

 For $s \in [0,S_1)$, the solution to equations \eqref{N-U}-\eqref{e-1-2} is given as $\gamma_0(s) =\tilde  \Lambda_0^{-1}(s)$, where
$ \tilde \Lambda_0(u) = \int_0^u f(v,0) dv$ for all $u > 0$; and, writing $N_{0,1} = N,\ Y_{0,1} =Y,$
$$N(t) = Y(t) = 0 \mbox{ for all } t \in [0, \gamma_0(S_1)).$$
Next, for $s \in  [S_1,S_2)$, $\gamma_0(s) = \tilde \Lambda_1^{-1}(s)$, where
\begin{equation*}
\tilde \Lambda_1(u)=\begin{cases}
\tilde \Lambda_0(u), \mbox{ if } u \in [0, \gamma_0(S_1))\\
S_1 + \int_{\gamma_0(S_1)}^u f(v,\tau_1\wedge(v-\gamma_0(S_1))) dv, \mbox{ if } u \ge \gamma_0(S_1),
\end{cases}
\end{equation*}
and
$$N(t) = 1, \;  Y(t) = \tau_1 \wedge (t-\gamma_0(S_1)) \mbox{ for all } t \in [\gamma_0(S_1), \gamma_0(S_2)).$$
One can now similarly write expressions for  $N(t)$ and $Y(t)$ for $t \in [\gamma_0(S_i), \gamma_0(S_{i+1}))$, for $i >1$.

This recursive construction shows in particular that $\gamma_0$ is a $\{\mathcal{F}_u\}$-adapted process, where 
\[
\mathcal{F}_u=\sigma\{ \xi_i(A): A\in \mathcal{B}([0,u]\times [0,\infty)), i=1,\dots,d\}.
\]
Consequently, for any $t \ge 0$, $\Lambda_0(t) = \gamma_0^{-1}(t)$ is a bounded $\{\mathcal{F}_u\}$-stopping time and therefore,
\[
N_{0,i}(t)-\Lambda_0(t)=\tilde{\xi}_{0,i}([0,\Lambda_0(t)]\times[0,\infty))
\]
is a $\{\mathcal{G}_t\}=\{\mathcal{F}_{\Lambda_0(t)}\}$-martingale, where $\tilde{\xi}_{0,i}=\xi_{0,i}-\eta$ is the compensated Poisson random measure associated with $\xi_{0,i}$, $i=1,\dots,d$.%
\subsection{Scaled Workload and Main Results}\label{mainres}
We now introduce the scaled system.
%
For each fixed $n \in \NN$, let $\xi_{n,1},\dots,\xi_{n,d}$ be independent Poisson random measures on $[0,\infty)\times [0,\infty)$ with common  intensity measure 
\begin{equation}
\eta_n(ds,dr)=n^{\alpha}ds \nu_n(dr),
\end{equation}
 where $\nu_n$ is introduced in \eqref{dist-2}.
%
Define, for $i=1,\dots,d$,
\begin{eqnarray}
X_{n,i}(t)&=& \frac{1}{n^\alpha}N_{n,i}(t) = \frac{1}{n^\alpha}\xi_{n,i}(B_n(t)),\nonumber\\
Y_{n,i}(t)&=&\frac{1}{n^{\alpha-1}}\int_{B_n(t)}r\wedge(t-\gamma_n(s))\xi_{n,i}(ds,dr)\label{e-2-2},
\end{eqnarray}
where
\begin{align}
\Lambda_n(t) =&\  \int_0^t f(s,\bar{Y}_n(s))ds,\  \bar{Y}_n(t) = \frac{1}{d} \sum_{i=1}^d Y_{n,i}(t),\  \gamma_n(t) = \  \Lambda_n^{-1}(t),\label{new-new-2-8}
 \end{align}
 and
\[
B_n(t)= [0, \Lambda_n(t)] \times [0, \infty).
\]

Note that $B_n$ and $\gamma_n$ do not depend on $i=1,\dots,d$.
As for the unscaled system, we  see that the solution $(\mathbf{X}_n,\mathbf{Y}_n)^T$ of the system (\ref{e-2-2}) exists and is unique on $[0,\infty)$ for each $n$, where $\mathbf{X}_n=(X_{n,1},\dots,X_{n,d})^T $ and $\mathbf{Y}_n=(Y_{n,1},\dots,Y_{n,d})^T$; and moreover, $\mathbf{X}_n$, $\mathbf{Y}_n\in D_{\mathbb{R}_+^d}[0, \infty)$.

Consider the ODE
\begin{equation}
\begin{cases}
\dot U(t) = a f(t,U(t)),\; t \ge 0, \; \\
U(0) = 0,\label{e-2-4}
\end{cases}
\end{equation}
where 
$$a= (\beta-1)\theta\int_0^\infty r(\theta r+1)^{-\beta}dr= (\beta-1)\theta\int_0^\infty n^2r(n\theta r+1)^{-\beta}dr=\frac{1}{\theta(\beta-2)}.$$ 
The following proposition will be proved in Section \ref{sec:sec3}.
\begin{proposition}\label{L-3-1}
 There is a unique continuous function $U$ that solves \eqref{e-2-4}.  The solution satisfies
 $$ \sup_{t \ge 0} |U(t)-bt|  < \infty.$$
 In the case when $b=a$, we have
$U(t) = bt$,  for all $t \ge 0.$
\end{proposition}

\begin{remark}\label{R-3-2}
As an immediate consequence of the above proposition we have that $f(t, U(t))=\exp\{-g(U(t)-bt)\}$ is bounded above and bounded below away from $0$, namely
\[
0<\inf_{t\geq0}\{f(t,U(t))\}\leq \sup_{t\geq 0}\{f(t,U(t))\}<\infty.
\]
 Denote
\[
f_y(t,y)=\frac{\partial}{\partial y}f(t,y)=-\exp\{-g(y-bt)\}g'(y-bt).
\]
Since $g'$ is bounded from below and above, it follows from the above proposition that $f_y(t,U(t))$ is also bounded below and bounded above away from $0$, namely
\[
-\infty<\inf_{t\geq0}\{f_y(t,U(t))\}\leq \sup_{t\geq 0}\{f_y(t,U(t))\}<0.
\]
\end{remark}
\bigskip

Let $\mathbf{X} = (X_1,\dots,X_d)^T=a^{-1}(U, U, \dots, U)^T$, $\mathbf{Y} =(Y_1,\dots,Y_d)^T= (U, U, \dots, U)^T$. The following is the first main result of this work.
\begin{theorem}\label{conv-prob}
As $n \to \infty$, $(\mathbf{X}_n, \mathbf{Y}_n)^T \to (\mathbf{X},\mathbf{Y})^T$ in $D_{\RR_+^{2d}}[0, \infty)$, in probability. Furthermore, for any $t>0$ and any $q \in [0, \beta-2)$,
 \begin{equation}\label{eq:strong117}
\sup_{0\le s \le t} n^{q} |\bar{Y}_n(s) - U(s)| \to 0,
\end{equation}
in probability, as $n\to\infty$. 
\end{theorem}

Let $V \in C_{\RR}[0,\infty)$ be given as the solution of
\begin{equation}\label{e-3-19-1}
V(t)=a\int_0^tf_y(s, U(s))V(s)ds-a\theta^{2-\beta}\int_0^t\frac{f(s,U(s))}{(t-s)^{\beta-2}}ds.
\end{equation}
From Remark \ref{R-3-2} the solution $V$ of the above linear equation exists and is unique. 

Define 
\begin{equation}
	\label{eq:ab348}
Z_{n,i}(t)=n^{\frac{\alpha+\beta-3}{2}}\left(Y_{n,i}(t)-Y_i(t)-\frac{V(t)}{n^{\beta-2}}\right), \; i=1,\dots, d.
\end{equation} 
Our next result gives the limiting behavior of the processes $\mathbf{Z}_n=(Z_{n,1},\dots,Z_{n,d})^T$. Note that $Y_i(t)+\frac{V(t)}{n^{\beta-2}}$ is not the expectation of $Y_{n,i}$, and hence, $Z_{n,i}$ in the above equation is not the conventional centered process of $Y_{n,i}$. However, from Proposition \ref{L-3-1} and Theorem \ref{conv-prob}, one can show $Y_i(t)=\lim\limits_{n\to\infty}\EE(Y_{n,i}(t))$. Also, as $n$ increases to infinity, the term $\frac{V(t)}{n^{\beta-2}}$ tends to zero. Thus, the next result can be regarded as a central limit theorem for the scaled and (nearly) centered process $\mathbf{Y}_n$.

\begin{theorem}\label{CLT}
 As $n \to \infty$, $\mathbf{Z}_n$ converges in distribution in $D_{\RR^d}[0, \infty)$ to $\mathbf{Z} = (Z_1, \dots, Z_d)^T$  where $\mathbf{Z}$ satisfies
\begin{equation}\label{linear}
 Z_i(t)=\int_{B(t)}r\wedge(t-\gamma(s))\Sigma_i(ds,dr)+a\int_0^tf_y(s,U(s))\bar{Z}(s)ds,\,i=1,\dots,d,
\end{equation}
where
\begin{equation}\label{new-new-2-14}
\Lambda(t) =    \int_0^t f(s,U(s))ds, \; \gamma(t) = \Lambda^{-1}(t),
\end{equation}
\[
B(t)=   [0, \Lambda(t)] \times [0, \infty), 
\]
\[
\bar{Z}(t) = \frac{1}{d}\sum_{i=1}^d Z_i(t),
\] and $\Sigma_1,\dots,\Sigma_d$ are independent Gaussian random measures on $[0,\infty)\times[0,\infty)$ with common control measure
$ds  (\beta-1)\theta^{1-\beta} r^{-\beta} dr$.
 \end{theorem} 

Integrals with respect to Gaussian random measures characterized by a control measure are defined, for example, in Chapter 3 of \cite{ST}.

\begin{remark}\label{Remark-2-2}
When $\lambda\equiv 1$ (constant) or $f\equiv 1$, note that $\Lambda_n(t)=t$ and $\gamma_n(t)=t$ in \eqref{new-new-2-8}, $B_n(t)=[0,t]\times[0,\infty)$ and hence
\begin{equation*}
Y_{n,i}(t)=\frac{1}{n^{\alpha-1}}\int_0^t\int_0^\infty r\wedge(t-s)\xi_{n,i}(ds,dr)
\end{equation*}
in \eqref{e-2-2}. After the change of variables $s\to s/n,\ r\to r/n$, this can be written as
\begin{align}\label{zeta}
 Y_{n,i}(t)=&\ \frac{1}{n^\alpha}\int_0^{nt}\int_0^\infty r\wedge(nt-s)\xi_{n,i}\left(d\left(\frac{s}{n}\right),d\left(\frac{r}{n}\right)\right)\nonumber\\
=&\ \frac{1}{n^\alpha}\int_0^{nt}\int_0^\infty r\wedge(nt-s)\zeta_{n,i}\left(ds,dr\right),
\end{align}
where $\zeta_{n,i}$ is a Poisson random measure with intensity measure $n^{\alpha-1}ds\ (\beta-1)\theta(\theta r+1)^{-\beta}dr$.  Written as \eqref{zeta}, $n^\alpha Y_{n,i}$ can be interpreted as the cumulative workload in the system scaled in time by $n$ and where heavy-tailed workloads are associated with sources arriving at Poisson rate $\lambda_n=n^{\alpha-1}$. This is the view taken, for example, in \cite{KT, MRRS}. It is well known that, after proper normalization and centering, the total workload converges to fractional Brownian motion in the so-called fast regime, that is, when
\[
\frac{\lambda_n}{n^{(\beta-1)-1}}=\frac{n^{\alpha-1}}{n^{\beta-2}}=n^{\alpha-\beta+1}\to\infty.
\]
This holds when $\alpha-\beta+1>0$, which is a part of our assumption \eqref{eq:ab1210}. It is also known that the normalization of the right-hand side of \eqref{zeta} (to the central limit theorem) is
\[
\frac{n^\alpha}{(\lambda_nn^{3-(\beta-1)})^{1/2}}=\frac{n^\alpha}{n^{(\alpha-\beta+3)/2}}=n^{\frac{\alpha+\beta-3}{2}},
\]
which coincides with that used in \eqref{eq:ab348}.
\end{remark}

\begin{remark} \label{R-2-2} Let $\mathbf{Z}^* = (Z_1^*, \dots, Z_d^*)^T$ be given as the solution of 
	\begin{equation}\label{linear-1}
	Z_i^*(t)=  R_i^*(t) + a\int_0^t f_y(s, U(s))\bar Z^*(s)ds,
	\end{equation}	
	where
	$$R_i^*(t)=\int_0^t\int_0^\infty \left( f(s, U(s))\right)^{1/2} (r\wedge(t-s))\Sigma_i(ds,dr), \; i=1,\dots, d,$$
	and \[
	\bar Z^*(t)=\frac{1}{d}\sum_{i=1}^dZ_i^*(t).
	\]
	One can check that $\mathbf{R}^* = (R_1^*, \dots, R_d^*)^T$  and $\mathbf{R} = (R_1, \dots, R_d)^T$ have the same distribution, where
	 $$R_i(t)=\int_{B(t)}r\wedge(t-\gamma(s))\Sigma_i(ds,dr), \; i=1,\dots, d.$$
	 Consequently, $\mathbf{Z}$ and $\mathbf{Z}^*$ are equal in law and thus \eqref{linear-1} gives an alternative representation for the weak limit of $\mathbf{Z}_n$, as $n\to \infty$.
	\end{remark}
	The following result shows the moment stabilization property of the admission control policy $g$.
\begin{theorem}
		\label{momstab}
The following uniform moment bound holds: 
	\begin{equation}
	\sup_{t \ge 0}\EE(\vert \bar Z(t)\vert^2)
	\leq \frac{2\theta^{1-\beta}}{d(\beta-2)(3-\beta)(a\mu)^{4-\beta}}\Gamma(4-\beta),
	\end{equation} 
	 where $\mu:=\inf_{s\geq0}\{-f_y(s,U(s))\} \in (0, \infty)$  and $\Gamma(\cdot)$ is the Gamma function.
	\end{theorem}
\begin{remark}
The case when there is no admission control corresponds to $g\equiv 0$.  Although the function $g=0$ does not satisfy Assumption \ref{ass:ass1}, it can be
shown along similar lines that in this case Theorem \ref{conv-prob} holds with $U(t) = at$, and therefore, $\sup_{t\ge0}|U(t) - bt|$ will be finite if and only
if $b=a$.  Furthermore, Theorem \ref{CLT} will hold as well (when $b=a$) but the moment stabilization property in Theorem \ref{momstab} fails.
\end{remark}
Finally we consider the asymptotic behavior of $\bar{Z}(T+ \cdot)$ as $T \to \infty$.  Here we restrict ourselves to  the case $b=a$. Then from Proposition \ref{L-3-1}, \eqref{linear} and \eqref{new-new-2-14}, the limit process in Theorem \ref{CLT} can be written as
\begin{equation}\label{new-2-11}
 Z_i(t)=\int_0^t\int_0^\infty r\wedge(t-s)\Sigma_i(ds,dr)-\kappa\int_0^t\bar{Z}(s)ds,\ i=1,\dots,d,
\end{equation}
where $\kappa=ag'(0)\in(0,\infty)$.

Let $B_{H}=(B_H(t),\, t\geq 0)$ be a standard fractional Brownian motion with Hurst parameter $H = \frac{4-\beta}{2}\in\left(\frac{1}{2},1\right)$, namely, $B_H$ is a mean zero Gaussian process with covariance 
\begin{equation}
\EE(B_{H}(t)B_{H}(s))=\frac{1}{2}\left(
t^{2H}+s^{2H}-|t-s|^{2H}\right) \,.  \label{cov}
\end{equation}%
 Let $Z_{\infty}(0)$ be a normal random variable with mean zero and variance 
\begin{equation}\label{var-Z-inf}
\sigma_0^2:=\EE(\vert Z_\infty(0)\vert^2)=\frac{\theta^{1-\beta}}{d\,(\beta-2)}\int_{0}^\infty\int_0^\infty e^{-\kappa v}e^{-\kappa u}\vert u-v\vert^{2-\beta}du\,dv<\infty,
\end{equation}
and let $(B_H, Z_{\infty}(0))$ be jointly Gaussian and 
the covariance function of $B_H$ and $Z_{\infty}(0)$ be 
\[
cov(B_H(t), Z_{\infty}(0))=\ \frac{\theta^{1-\beta}}{\sigma d\,(\beta-2)}\int_0^t\int_0^\infty e^{-\kappa v}(u+v)^{2-\beta}dv\,du,
\]
where  $\sigma=\sqrt{\frac{2\theta^{1-\beta}}{d\,(\beta-2)(3-\beta)(4-\beta)}}$.
 Let $Z_\infty$ be
the   fractional Ornstein-Uhlenbeck process given as the unique solution of
\begin{equation}\label{Z-2-16}
Z_\infty(t)=Z_\infty(0)-\kappa\int_0^tZ_\infty(s)ds+\sigma B_{H}(t).
\end{equation}

\begin{theorem} \label{ouplim}
 Let $b=a$ and let $\mathbf{Z}$ be as in Theorem \ref{CLT}. Then, as $T\to\infty$, $\bar Z(T+\cdot)$ converges in distribution, in $C_{\RR}[0,\infty)$
to $Z_{\infty}$ given by \eqref{Z-2-16}. Moreover, the process $Z_\infty$ is stationary.
\end{theorem}
\bigskip

\setcounter{equation}{0}
\section{Law of Large Numbers}  \label{sec:sec3}
In this section we will prove Proposition \ref{L-3-1} and Theorem \ref{conv-prob}.\\ \ 

\noindent {\bf Proof of Proposition \ref{L-3-1}:} 
Consider the ODE
\begin{equation}
	\label{eq:ab259}
	     \begin{cases}
	     \dot u(t) = a \exp\{-g(u(t))\} -b, \; t \ge 0, \;\\
	      u(0) = 0.
	      \end{cases}
\end{equation}
	Clearly, a differentiable function $u$ solves \eqref{eq:ab259} if and only if
	$U(t) = u(t)+ bt$ solves \eqref{e-2-4}.
	From Assumption \ref{ass:ass1}, the function $h(x) = a\exp\{-g(x)\} - b$, $x \in \RR$, is locally
	Lipschitz.  For each $n \in \NN$, define $h_n(x) = h((x\wedge n)\vee(-n))$, $x \in \RR$.  Since $h_n$ is a Lipschitz function, for any $n \in \NN$, the ODE
	\begin{equation}
	\label{eq:ab305}
	        \begin{cases}
	        \dot u(t) = h_n(u(t)), \; t \ge 0, \;\\
	         u(0) = 0
	         \end{cases}
	\end{equation}
	has a unique solution $u_n$.  Let $K$ be the unique solution of the equation
	$$a\exp(-g(K)) - b = 0, $$
	i.e.
	$g(K      )=\log{\frac{a}{b}}.$
	
Then, for all $n > |K|$, if $b > a$, $u_n(t) \le 0$ for all $t$ and $u_n(t)$ decreases to $K \in (-\infty , 0)$; 
if $b < a$, $u_n(t) \ge 0$ for all $t$ and $u_n(t)$ increases to $K \in (0 , \infty)$; and finally if $b=a$, $u_n(t)=0$ for all $t$.
Consequently, for any $n > |K|$, 
\begin{equation}
		\label{eq:unifbdd}
		\sup_{t\ge 0} |u_n(t)| \le |K|
			\end{equation}
 and $u_n$ solves \eqref{eq:ab259}.  This proves the existence of solutions.  

Now consider uniqueness.  Let $\tilde u$ be another solution of
\eqref{eq:ab259}.  Let $\tau = \inf\{t: |\tilde u(t)| \ge |K|+1\}$.  From unique solvability of \eqref{eq:ab305}, for any $n \ge |K|+1$,
$\tilde u(t) = u_n(t)$ for all $t \in [0, \tau)$.  From \eqref{eq:unifbdd} we now see that $\tau = \infty$.  This proves unique solvability
of \eqref{eq:ab259} and consequently that of \eqref{e-2-4}.  Also, as noted above,
$$\sup_{t\ge 0} |U(t)-bt| = \sup_{t\ge 0} |u(t)| \le |K |$$
and
$
U(t) - bt = u(t) = 0$ for all $t$, if $b=a$.  The result follows. $\Box$
\\ \

Next, we give the proof of Theorem \ref{conv-prob}.\\ \ 

\noindent {\bf Proof of Theorem \ref{conv-prob}:}
Let $\tilde{\xi}_{n,i}=\xi_{n,i}-\eta_n$ be the compensated Poisson random measure associated with $\xi_{n,i}$, $i=1,\dots,d$.
Rewrite $\mathbf{X}_n, \mathbf{Y}_n$ as 
\begin{eqnarray}
X_{n,i}(t)&=&\frac{1}{n^\alpha}\tilde{\xi}_{n,i}(B_n(t))+\frac{1}{n^\alpha}\eta_n(B_n(t))=\frac{1}{n^\alpha}\tilde{\xi}_{n,i}(B_n(t))+\Lambda_n(t),\nonumber\\
Y_{n,i}(t)&=&\frac{1}{n^{\alpha-1}}\int_{B_n(t)}r\wedge(t-\gamma_n(s))\tilde{\xi}_{n,i}(ds,dr)\nonumber\\
&&+\frac{1}{n^{\alpha-1}}\int_{B_n(t)}r\wedge(t-\gamma_n(s))n^\alpha (\beta-1)n\theta(n\theta r+1)^{-\beta}dsdr.\label{eq:ab1703}
\end{eqnarray}
By the change of variables $s=\int_0^vf(u,\bar{Y}_n(u))du=\Lambda_n(v)$, the second term on the right-hand side of \eqref{eq:ab1703} equals
\begin{equation}\label{s-s-2-8}
n^2\theta (\beta-1)\int_0^t\int_0^\infty f(v, \bar{Y}_n(v))\left(r\wedge (t-v)\right)(n\theta r+1)^{-\beta}drdv.
\end{equation}
Consider the inner integral in \eqref{s-s-2-8}. For $0\leq v<t$, by changing variables, we see that
\begin{eqnarray}
&&\theta (\beta-1)n^2 \int_0^\infty\left(r\wedge (t-s)\right)(n\theta r+1)^{-\beta}dr\nonumber\\
&=&\theta (\beta-1)n^2\left [\int_0^{t-s}r(n\theta r+1)^{-\beta}dr+\int_{t-s}^\infty (t-s)(n\theta r+1)^{-\beta}dr\right]\nonumber\\
&=&\theta (\beta-1)\left[\int_0^{n(t-s)}r(\theta r+1)^{-\beta}dr+n(t-s)\int_{n(t-s)}^\infty (\theta r+1)^{-\beta}dr\right]\nonumber\\
&=&\frac{1}{\theta (\beta-2)}-\frac{1}{\theta(\beta-2)(n\theta(t-s)+1)^{\beta-2}}=a\left(1-\frac{1}{(n\theta(t-s)+1)^{\beta-2}}\right).
\end{eqnarray}
Therefore, for each $i=1,\dots,d$,
\begin{eqnarray}
Y_{n,i}(t)&=&\frac{1}{n^{\alpha-1}}\int_{B_n(t)}r\wedge(t-\gamma_n(s))\tilde{\xi}_{n,i}(ds,dr)\nonumber\\
&&+\int_0^taf(s, \bar{Y}_n(s))\left(1-\frac{1}{(n\theta(t-s)+1)^{\beta-2}}\right)ds
\end{eqnarray}
and hence,
\begin{eqnarray}
X_{n,i}(t)-X_i(t) &= &\frac{1}{n^\alpha}\tilde{\xi}_{n,i}(B_n(t))+\int_0^tf(s,\bar{Y}_n(s))ds-\int_0^tf(s,U(s))ds,\label{eq:ab1718}\\
Y_{n,i}(t)-Y_i(t)&=& \frac{1}{n^{\alpha-1}}\int_{B_n(t)}r\wedge(t-\gamma_n(s))\tilde{\xi}_{n,i}(ds,dr)\nonumber\\
&&+\int_0^t a[f(s,\bar{Y}_n(s))-f(s, U(s))]ds-a\int_0^t\frac{f(s,\bar{Y}_n(s))}{(n\theta(t-s)+1)^{\beta-2}} ds.\nonumber\\ \label{e-2-6}
\end{eqnarray}

We will first show \eqref{eq:strong117}.
Since $\bar{Y}_n(t)\geq 0$ for all $t\geq0$ and $\beta \in (2,3)$, from \eqref{new-eq-2-4}, we have that
\begin{equation*}
\int_0^t\frac{f(s,\bar{Y}_n(s))}{(n\theta (t-s)+1)^{\beta-2}} ds\leq \frac{1}{n^{\beta-2}\theta^{\beta -2}}\int_0^t\frac{\exp\{-g(-bs)\}}{(t-s)^{\beta-2}}ds\leq 
\frac{\left(\sup_{0\leq s\leq t}\exp\{-g(-bs)\}\right)t^{3-\beta}}{n^{\beta-2}\theta^{\beta-2}(3-\beta)}.
\end{equation*}
Consequently, for every $t> 0$, 
\begin{equation}
\label{e-3-10}	
\lim_{n\to \infty} n^{q}\int_0^t\frac{f(s,\bar{Y}_n(s))}{(n\theta (t-s)+1)^{\beta-2}} ds = 0, \mbox{ a.s.}, \mbox{ for all } q \in [0, \beta -2).
\end{equation}

For $n \in \NN$, define the filtration $\{\mathcal{F}_u^n\}$ as
\[
\mathcal{F}^n_u=\sigma\{\xi^i_n(A):\,A\in\mathcal{B}([0,u]\times[0,\infty)),i=1,\dots,d\}.
\]
Then for each $i=1,\dots,d$, $\tilde{\xi}_{n,i}([0,u]\times[0,\infty))$ is an $\{\mathcal{F}^n_u\}$-martingale.
As for the unscaled process in Section \ref{sec:sec2.1}, $\gamma_n$ is a continuous, strictly increasing $\{\mathcal{F}^n_u\}$-adapted process.
Consequently, for every $t \ge 0$, $\Lambda_n(t)=\gamma_n^{-1}(t)$ is a $\{\mathcal{F}^n_u\}$-stopping time.
Now consider the first term on the right-hand side of \eqref{e-2-6}, that is, for  $i=1,\dots,d$,
\[
 \frac{1}{n^{\alpha-1}}\int_{B_n(t)}r\wedge(t-\gamma_n(s))\tilde{\xi}_{n,i}(ds,dr)=:A_{n,i}(t).
\]
Observe that 
\[
A_{n,i}(t)=\clu_{n,i}^{(1)}(\Lambda_n(t)),
\]
where, for $i=1,\dots,d$,
$$\clu_{n,i}^{(1)}(u) = \frac{1}{n^{\alpha-1}}\int_{[0,u]\times [0, \infty)} r\wedge(t-\gamma_n(s))_+\tilde{\xi}_{n,i}(ds,dr).$$
Note that  $\clu_{n,i}^{(1)}(u)$ is a $\{\mathcal{F}^n_u\}$-martingale with predictable quadratic variation process
$$\frac{1}{n^{2(\alpha-1)}}\int_{[0,u]\times [0, \infty)} \left (r\wedge(t-\gamma_n(s))_+\right)^2 \eta_{n}(ds,dr).$$
Using a change of variables, we have, for each $t>0$ and $\mathbf{A}_n=(A_{n,1},\dots,A_{n,d})^T$,
\begin{align}\label{e-2-9}
&\EE(\vert \mathbf{A}_n(t)\vert^2)\nonumber\\
=&\ \EE\left(\sum_{i=1}^d\left\vert \clu_{n,i}^{(1)}(\Lambda_n(t))\right\vert^2\right)\nonumber\\
=&\ \frac{1}{n^{2\alpha-2}}\EE\left(\sum_{i=1}^d\int_{[0,\Lambda_n(t)]\times[0,\infty)} 
\left(r\wedge(t-\gamma_n(s))\right)^2\eta_n(ds,dr)\right)\nonumber\\
=&\ \frac{n\theta (\beta-1)d}{n^{\alpha-2}}\EE\left(\int_0^t\int_0^\infty f(s,\bar{Y}_n(s))(r\wedge(t-s))^2(n\theta r+1)^{-\beta}dr\,ds\right)\nonumber\\
\leq&\ \frac{n\theta (\beta-1)d}{n^{\alpha-2}}\int_0^t\int_0^{t-s}\exp\{-g(-bs)\}r^2(n\theta r+1)^{-\beta}dr\,ds\nonumber\\
&\ +\frac{n\theta (\beta-1)d}{n^{\alpha-2}}\int_0^t\int_{t-s}^\infty \exp\{-g(-bs)\}(t-s)^2(n\theta r+1)^{-\beta}dr\,ds\nonumber\\
\leq&\ \frac{\theta (\beta-1)d\ \sup_{0\leq s\leq t}\{\exp\{-g(-bs)\}\}}{n^\alpha}\int_0^t\int_0^{n(t-s)}r^2(\theta r+1)^{-\beta}dr\,ds\nonumber\\
&\ +\frac{\theta (\beta-1)d\ \sup_{0\leq s\leq t}\{\exp\{-g(-bs)\}\}}{n^\alpha}\int_0^t\int_{n(t-s)}^\infty (n(t-s))^2(\theta r+1)^{-\beta}dr\,ds.\nonumber\\
\end{align}
For the first term on the right-hand side of \eqref{e-2-9}, note that
\begin{eqnarray*}\label{s-e-2-9}
\frac{1}{n^\alpha}\int_0^t\int_0^{n(t-s)}r^2(\theta r+1)^{-\beta}dr\,ds&\leq&\frac{1}{\theta^{\beta}n^\alpha}\int_0^t\int_0^{n(t-s)}r^{2-\beta}dr\,ds\nonumber\\
&=&\frac{1}{\theta^{\beta}(3-\beta)}n^{3-\beta-\alpha}\int_0^t(t-s)^{3-\beta}ds\nonumber\\
&=&\frac{n^{3-\beta-\alpha}t^{4-\beta}}{\theta^{\beta}(3-\beta)(4-\beta)}.
\end{eqnarray*}
Since $\beta -2 < \frac{\beta + \alpha -3}{2}$ (or equivalently, $\alpha>\beta-1$), we obtain that, for all $t \ge 0$,
\begin{equation}
	\label{eq:cgcet1}
	\frac{n^{2q}}{n^\alpha}\int_0^t\int_0^{n(t-s)}r^2(\theta r+1)^{-\beta}dr\,ds \to 0 , \mbox{ for all } q \in [0, \beta-2).
\end{equation}
For the second term on the right-hand side of  \eqref{e-2-9}, note that
\begin{eqnarray*}\label{s-o-e-2-9}
&&\frac{1}{n^\alpha}\int_0^t\int_{n(t-s)}^\infty(n(t-s))^2(\theta r+1)^{-\beta}dr\,ds\nonumber \\
&\leq& 
\frac{1}{\theta^{\beta}n^\alpha}\int_0^t\int_{n(t-s)}^\infty(n(t-s))^2r^{-\beta}dr\,ds\nonumber\\
&=&\frac{n^{3-\beta-\alpha}}{\theta^{\beta}(\beta-1)}\int_0^t(t-s)^{3-\beta}ds=\frac{n^{3-\beta-\alpha}t^{4-\beta}}{\theta^{\beta}(\beta-1)(4-\beta)}.
\end{eqnarray*}
Thus we have
\begin{equation}
	\label{eq:cgcet2}
	\frac{n^{2q}}{n^\alpha}\int_0^t\int_{n(t-s)}^\infty(n(t-s))^2(\theta r+1)^{-\beta}drds \to 0 , \mbox{ for all } q \in [0, \beta-2).
\end{equation}
 Combining \eqref{e-2-9}, \eqref{eq:cgcet1} and \eqref{eq:cgcet2} we conclude that
\begin{equation}\label{s-2-17}
\lim_{n\to \infty}n^{2q}\EE\left(\vert \mathbf{A}_n(t)\vert^2\right) = 0,\ \mbox{ for all } q \in [0, \beta-2).
\end{equation}

We argue next that $n^{q}\mathbf{A}_n=n^q(A_{n,1},\dots,A_{n,d})^T$ converges to the zero process, in $D_{\RR^d}[0,\infty)$, in probability. In view of \eqref{s-2-17},
it suffices to check that $\{n^{q}\mathbf{A}_n\}$ is tight.  To prove tightness we will use a standard tightness criterion.
Namely, we will show that  for each fixed $T>0$ there exists $C_T >0$ such that for $0\leq h\leq 1$ and $h\leq t\leq T$,
\begin{equation}\label{s-s-2-22}
n^{4q}\EE\left(\vert \mathbf{A}_n(t+h)-\mathbf{A}_n(t)\vert^2\vert\mathbf{A}_n(t)-\mathbf{A}_n(t-h)\vert^2\right)\leq C_Th^2.
\end{equation}
The above inequality, together with the relative compactness of $n^{q}\mathbf{A}_n(t)$ for each $t\geq0$ (which follows from \eqref{s-2-17}),
yields tightness of $\{n^{q}\mathbf{A}_n\}$ (cf. Theorems 3.8.6 and 3.8.8 in  \cite{EK}).

Now fix $T>0$. In order to show (\ref{s-s-2-22}), it is sufficient to prove that, for any $0\leq h\leq 1$ and $0\leq t\leq T$,
\begin{equation}\label{eq:ab1933}
n^{4q}\EE\left(\vert \mathbf{A}_n(t+h)-\mathbf{A}_n(t)\vert^4\right)\leq C_Th^2.
\end{equation}
In the following, we use $C_T>0$ to denote a generic constant depending on $T, \ \theta$ and $\beta$ whose value may vary from line to line. 
For $r,s,h,t \in \RR$, denote
\begin{equation}\vartheta_n^{h,t}(r,s) =  r\wedge(t+h-\gamma_n(s))_{+}-r\wedge(t-\gamma_n(s))_{+}. \label{eq:vnhtrs}
	\end{equation}

Define, for $i=1,\dots,d$,
\begin{equation}\label{eq:clu3}
\clu^{(2)}_{n,i}(u)=\frac{1}{n^{\alpha-1}}\int_{[0,u]\times[0,\infty)}\vartheta_n^{h,t}(r,s)\tilde{\xi}_{n,i}(ds,dr).
\end{equation}
Observe that  $\clu^{(2)}_{n,i}(u)$ is a $\{\mathcal{F}^n_u\}$-martingale  with quadratic variation process
$$\frac{1}{n^{2(\alpha-1)}}\int_{[0,u]\times [0, \infty)}\left(\vartheta_n^{h,t}(r,s)\right)^2 \xi_{n}(ds,dr).$$
Since  $\gamma_n(s)\leq t$ if and only if $s\leq \Lambda_n(t)$, we have
 $$A_{n,i}(t+h)-A_{n,i}(t)=\clu^{(2)}_{n,i}(\Lambda_n(t+h)).$$
 Recalling that $\Lambda_n(t+h)$ is a $\{\mathcal{F}^n_u\}$-stopping time, we have by the Burkholder-Davis-Gundy inequality 
that for some $C>0$, 
\begin{eqnarray}\label{s-2-21}
&&\EE\left(\vert A_{n,i}(t+h)-A_{n,i}(t)\vert^4\right)=\EE\left(\vert \clu^{(2)}_{n,i}(\Lambda_n(t+h))\vert^4\right)\nonumber\\
&\leq& \frac{C}{n^{4(\alpha-1)}} \EE\left(\left(\int_{[0,\Lambda_n(t+h)]\times[0,\infty)}[\vartheta_n^{h,t}(r,s)]^2\xi_{n,i}(ds,dr)\right)^{2}\right)\nonumber\\
&\leq&\frac{2C}{n^{4(\alpha-1)}} \EE\left(\left(\int_{B_n(t+h)}[\vartheta_n^{h,t}(r,s)]^2\tilde{\xi}_{n,i}(ds,dr)\right)^{2}\right)\nonumber\\
&&+\frac{2C}{n^{4(\alpha-1)}}\EE\left(\left(\int_{0}^{\Lambda_n(t+h)}\int_0^\infty[\vartheta_n^{h,t}(r,s)]^2n^\alpha\theta (\beta-1)n(n\theta r+1)^{-\beta}dr\,ds\right)^{2}\right)\nonumber\\
&=&\frac{2C}{n^{4(\alpha-1)}} \EE\left(\int_{0}^{\Lambda_n(t+h)}\int_0^\infty[\vartheta_n^{h,t}(r,s)]^4n^\alpha\theta (\beta-1)n(n\theta r+1)^{-\beta}dr\,ds\right)\nonumber\\
&&+\frac{2C}{n^{4(\alpha-1)}}\EE\left(\left(\int_{0}^{\Lambda_n(t+h)}\int_0^\infty[\vartheta_n^{h,t}(r,s)]^2n^\alpha\theta (\beta-1)n(n\theta r+1)^{-\beta}dr\,ds\right)^{2}\right).\nonumber\\
\end{eqnarray}
Denote for $r,s,h,t \in \RR$,
$$\tilde \vartheta_n^{h,t}(r,s) =  r\wedge(t+h-s)_{+}-r\wedge(t-s)_{+}.$$
By a change of variables and using (\ref{s-4-1}) from Lemma \ref{L-A-1} in Appendix we can  bound the first term on the right-hand side of
(\ref{s-2-21}) as
\begin{eqnarray}\label{s-2-22}
&&\frac{2}{n^{4(\alpha-1)}} \EE\left(\int_{0}^{\Lambda_n(t+h)}\int_0^\infty[\vartheta_n^{h,t}(r,s)]^4n^\alpha\theta (\beta-1)n(n\theta r+1)^{-\beta}dr\,ds\right)\nonumber\\
&=&\frac{2}{n^{4(\alpha-1)}} \EE\left(\int_{0}^{t+h}\int_0^\infty  f(s,\bar{Y}_n(s))[\tilde \vartheta_n^{h,t}(r,s)]^4n^\alpha\theta (\beta-1)n(n\theta r+1)^{-\beta}dr\,ds\right)\nonumber\\
&\leq&\frac{2}{n^{4(\alpha-1)}} \int_{0}^{t+h}\int_0^\infty  \exp\{-g(-bs)\}[\tilde \vartheta_n^{h,t}(r,s)]^4n^\alpha\theta (\beta-1)n(n\theta r+1)^{-\beta}dr\,ds\nonumber\\
&\leq&\frac{2\sup_{0\leq s\leq T+1}\{\exp\{-g(-bs)\}\}}{n^{4(\alpha-1)}} \int_{0}^{t+h}\int_0^\infty  [\tilde \vartheta_n^{h,t}(r,s)]^4n^\alpha\theta (\beta-1)n(n\theta r+1)^{-\beta}dr\,ds\nonumber\\
&\leq&C_Tn^{5-\beta-3\alpha}h^{6-\beta}.
\end{eqnarray}
For the second term on the right-hand side of (\ref{s-2-21}), by a change of variables once more and using (\ref{s-4-2}) in Lemma \ref{L-A-1}, we have
\begin{eqnarray}\label{s-2-23}
&&\frac{2}{n^{4(\alpha-1)}}\EE\left(\left(\int_{0}^{\Lambda_n(t+h)}\int_0^\infty[\vartheta_n^{h,t}(r,s)]^2n^\alpha\theta (\beta-1)n(n\theta r+1)^{-\beta}dr\,ds\right)^{2}\right)\nonumber\\
&=&\frac{2}{n^{4(\alpha-1)}}\EE\left(\left(\int_{0}^{t+h}\int_0^\infty f(s,\bar{Y}_n(s))[\tilde \vartheta_n^{h,t}(r,s)]^2n^\alpha\theta (\beta-1)n(n\theta r+1)^{-\beta}dr\,ds\right)^{2}\right)\nonumber\\
&\leq&\frac{2\sup_{0\leq s\leq T+1}\{e^{-2g(-bs)}\}}{n^{4(\alpha-1)}}\left(\int_{0}^{t+h}\int_0^\infty [\tilde \vartheta_n^{h,t}(r,s)]^2n^\alpha\theta (\beta-1)n(n\theta r+1)^{-\beta}dr\,ds\right)^{2}\nonumber\\
&\leq& C_Tn^{6-2\beta-2\alpha}h^{2(4-\beta)}.
\end{eqnarray}
Observing that
$$
\beta -2 < \frac{1}{4} \min\{\beta + 3\alpha-5, 2\beta + 2\alpha -6\} \mbox{ and } \min \{ 6-\beta , 2(4-\beta) \} > 2,$$
and combining \eqref{s-2-22} and \eqref{s-2-23}, we conclude that
\eqref{eq:ab1933} holds for every $q \in [0, \beta -2)$.
This shows that for every such $q$, $n^{q} \mathbf{A}_n$ converges in probability to the zero process, in $D_{\RR^d}[0,\infty)$ and thus for every $t > 0$
\begin{equation}\label{e-2-13}
\sup_{0\leq s\leq t} n^{q}\vert \mathbf{A}_n(s)\vert \to 0, \mbox{ as } n \to \infty , \mbox{ in probability, for every } q \in [0, \beta -2).
\end{equation}

Combining  \eqref{e-3-10}	 and \eqref{e-2-13} and recalling that $\frac{1}{d} \sum_{i=1}^d (Y_{n,i} - Y_i) = \bar{Y}_n-U$, we have from \eqref{e-2-6}
that
\begin{equation}\label{eq:ab1940}
	\bar{Y}_n(t) - U(t) = a \int_0^t \left [ f(s, \bar{Y}_n(s)) - f(s, U(s))\right] ds + \cls_n(t),
\end{equation}
where $$\cls_n(t)=\frac{1}{d}\sum_{i=1}^dA_{n,i}(t)-a\int_0^t\frac{f(s,\bar{Y}_n(s))}{(n\theta(t-s)+1)^{\beta-2}}ds=\bar{A}_n(t)-a\int_0^t\frac{f(s,\bar{Y}_n(s))}{(n\theta(t-s)+1)^{\beta-2}}ds$$ 
and $n^{q}\cls_n$ converges in probability to the zero process, in $D_{\RR^d}[0,\infty)$.
From Assumption \ref{ass:ass1}, we have that $y \mapsto f(t,y)$ is a Lipschitz function on $\RR_+$, uniformly in $t$ on compact intervals, since
\begin{equation}\label{new-3-25}
\sup_{y \in \RR_+}\vert f_y(t,y)\vert=\sup_{y \in \RR_+}\vert-\exp\{-g(y-bt)\}g'(y-bt)\vert\leq L\exp\{-g(-bt)\}
\end{equation}
 for all  $t\in[0,\infty)$.
Thus \eqref{eq:strong117} follows by an application of Gronwall's lemma to \eqref{eq:ab1940}.
%
%
%
%

Recall that for each $i=1,\dots,d$, $\tilde{\xi}_{n,i}([0,u]\times[0,\infty))$ is an $\{\mathcal{F}^n_u\}$-martingale and
 for every $t \ge 0$, $\Lambda_n(t)=\gamma_n^{-1}(t)$ is a $\{\mathcal{F}^n_u\}$-stopping time. Therefore,
 \begin{equation}\label{eq:ab1809}
M^{(1)}_{n,i}(t) =\tilde{\xi}_{n,i}([0,\Lambda_n(t)]\times[0,\infty))=\tilde{\xi}_{n,i}(B_n(t))
\end{equation}
is a $\{\mathcal{G}_t^n\}$-martingale, where  $\mathcal{G}_t^n=\mathcal{F}^n_{\Lambda_n(t)}$. By Doob's maximal inequality, for some $C>0$,
\begin{eqnarray}\label{e-2-8}
\PP\left(\sup_{0\leq s\leq t}\frac{1}{n^\alpha}\vert M^{(1)}_{n,i}(t)\vert\geq \epsilon \right)&\leq& \frac{C\EE\vert M^{(1)}_{n,i}(t)\vert^2}{n^{2\alpha} \epsilon^2}\nonumber\\
&=&\frac{C\EE\left(\eta_{n}([0, \Lambda_n(t)] \times [0, \infty))\right)}{n^{2\alpha} \epsilon^2} \nonumber\\
&=& \frac{C\EE\left(\int_0^tf(s,\bar{Y}_n(s))ds\right)}{n^{\alpha} \epsilon^2}\nonumber\\
&\leq&  \frac{C\int_0^t\exp\{-g(-bs)\}ds}{n^{\alpha} \epsilon^2}.
\end{eqnarray}
Combining (\ref{e-2-8})  and \eqref{eq:ab1809} we have that, as $n\to \infty$,
\begin{equation}\label{eq:ab1938}
\PP\left(\sup_{0\leq s\leq t}\frac{1}{n^\alpha}\left\vert \tilde{\xi}_{n,i}(B_n(t))\right\vert\geq \epsilon \right)\to 0.
\end{equation}
Thus the first term on the right-hand side of \eqref{eq:ab1718} converges to the zero process, uniformly on compacts, in probability, as $n \to \infty$. Finally,  combining \eqref{new-3-25}, \eqref{eq:ab1938}, \eqref{e-3-10}	 and \eqref{e-2-13} (with $q=0$),  we deduce applying Gronwall's lemma to \eqref{eq:ab1718} and \eqref{e-2-6} that $(\mathbf{X}_n, \mathbf{Y}_n)^T \to (\mathbf{X},\mathbf{Y})^T$ in $D_{\RR_+^{2d}}[0, \infty)$, in probability, as $n \to \infty$.
$\Box$

\bigskip

\setcounter{equation}{0}
\section{Central Limit Theorem}\label{sec:sec4}
In this section we prove Theorem \ref{CLT}.
From (\ref{e-3-19-1}) and (\ref{e-2-6}), we can write  \eqref{eq:ab348} as
\begin{eqnarray}\label{e-3-19}
Z_{n,i}(t)&=&\int_{B_n(t)}(r\wedge(t-\gamma_n(s)))\Sigma_{n,i}(ds,dr)\nonumber\\
&&+a\int_0^tn^{\frac{\alpha+\beta-3}{2}} \left[f(s, \bar{Y}_n(s))-f(s,U(s))-f_y(s,U(s))\frac{V(s)}{n^{\beta-2}}\right]ds\nonumber\\
&&-n^{\frac{\alpha+\beta-3}{2}}\int_0^t\frac{a f(s,\bar{Y}_n(s))\ ds}{(n\theta (t-s)+1)^{\beta-2}}+n^{\frac{\alpha-\beta+1}{2}}\int_0^t\frac{a\theta^{2-\beta}f(s,U(s))\ ds}{(t-s)^{\beta-2}},
\end{eqnarray}
where, with $\sigma_n = n^{(\alpha-\beta+1)/2}$,
$$\Sigma_{n,i}(A)=\frac{n^{\frac{\alpha+\beta-3}{2}}}{n^{\alpha-1}}\tilde{\xi}_{n,i}(A)=\sigma_n^{-1}\tilde{\xi}_{n,i}(A)$$
 is a random signed measure on $[0,\infty)\times[0,\infty)$, $i=1,\dots,d$.  Note that 
$$Var(\Sigma_{n,i}(A))= n^{\beta-1} m\times \nu_n(A), \ i=1,\dots,d,\ \mbox{\text{for}}\ A \in \clb(\RR_+^2)\ \mbox{\text{with}}\ m\times \nu_n(A)<\infty.$$

Note also that
\begin{eqnarray}\label{new-4-27}
&&f(s, \bar{Y}_n(s))-f(s,U(s))-f_y(s,U(s))\frac{V(s)}{n^{\beta-2}}\nonumber\\
&=&(\bar{Y}_n(s)-U(s))\int_0^1[f_y(s,U(s)+x(\bar{Y}_n(s)-U(s)))-f_y(s,U(s))]dx\nonumber\\
&&+\left(\bar{Y}_n(s)-U(s)-\frac{V(s)}{n^{\beta-2}}\right)f_y(s,U(s)).
\end{eqnarray}
Thus the middle term on the right-hand side of \eqref{e-3-19} equals
\begin{eqnarray*}\label{new-4-28}
&&a\int_0^tn^{\frac{\alpha+\beta-3}{2}}(\bar{Y}_n(s)-U(s))\int_0^1[f_y(s,U(s)+x(\bar{Y}_n(s)-U(s)))-f_y(s,U(s))]dx ds\nonumber\\
&+&a\int_0^tf_y(s,U(s))\bar{Z}_n(s)ds,
\end{eqnarray*}
where, recall, $\bar{Z}_n(s)=\frac{1}{d}\sum_{i=1}^dZ_{n,i}(s)$.
Let
\begin{align}
R_{n,i}(t)&=\int_{B_n(t)}(r\wedge(t-\gamma_n(s)))\Sigma_{n,i}(ds,dr),\ i=1,\dots,d,\label{R-n}\\
C_n(t)&=a\int_0^tn^{\frac{\alpha+\beta-3}{2}}(\bar{Y}_n(s)-U(s))\int_0^1[f_y(s,U(s)+x(\bar{Y}_n(s)-U(s)))-f_y(s,U(s))]dx\,ds,\nonumber\\
 \label{C-n}\\
D_n(t)&=n^{\frac{\alpha-\beta+1}{2}}\int_0^t\frac{a\theta^{2-\beta}f(s,U(s))\ ds}{(t-s)^{\beta-2}}-n^{\frac{\alpha+\beta-3}{2}}\int_0^t\frac{a f(s,\bar{Y}_n(s))\ ds}{(n\theta(t-s)+1)^{\beta-2}}.
\label{D-n}
 \end{align}
Letting $\mathbf{\mathcal{R}}_n(t)=(\mathcal{R}_{n,1}(t),\dots, \mathcal{R}_{n,d}(t))^T$, where $\mathcal{R}_{n,i}(t)=R_{n,i}(t)+C_n(t)+D_n(t)$,  we can rewrite equation 
(\ref{e-3-19}) as 
\begin{equation}\label{new-4-30}
Z_{n,i}(t)=\mathcal{R}_{n,i}(t)+a\int_0^tf_y(s,U(s))\bar{Z}_n(s)ds,\ i=1,\dots,d.
\end{equation}
\bigskip

\noindent 		{\bf Proof of Theorem \ref{CLT}:}
Define $\psi: D_{\mathbb{R}^d}[0,\infty)\to D_{\mathbb{R}^d}[0,\infty)$ by 
\[
[\psi(x)]_i(t)=x_i(t)+a\int_0^tf_y(s,U(s))\overline{\psi(x)}(s)ds, \ i=1,\dots,d, \; x\in D_{\mathbb{R}^d}[0,\infty),
\]
 where $\overline{\psi(x)}=\frac{1}{d}\sum_{i=1}^d[\psi(x)]_i$. Then,  $\psi$ is a continuous mapping from $D_{\mathbb{R}^d}[0,\infty)$ to $D_{\mathbb{R}^d}[0,\infty)$.
Also, from \eqref{new-4-30} we see that $\mathbf{Z}_n=(Z_{n,1},\dots, Z_{n,d})^T = \psi(\mathcal{R}_n)$.

  Combining Lemmas \ref{anto0}, \ref{bnto0} and \ref{sicgce} below, we see that
$\mathcal{R}_n$ converges to $\mathbf{R}=(R_1,\dots,R_d)^T$, in distribution, in $D_{\mathbb{R}^d}[0,\infty)$, where
\begin{equation}\label{R}
R_i(t) = \int_{B(t)}(r\wedge(t-\gamma(s)))\Sigma_{i}(ds,dr), \; i = 1, \dots, d,
\end{equation}
and $\Sigma_i$, $i=1,\dots,d$, is as in Theorem \ref{CLT}. The result now follows from continuous mapping theorem. $\Box$\\ \ 

The next three lemmas were used in the proof of Theorem \ref{CLT} above.
\begin{lemma}\label{anto0}
Let $C_n$ be as given in \eqref{C-n}. 	As $n \to \infty$, $\sup_{0 \le s \le t} |C_n(s)| \to 0$, in probability, for every $t \ge 0$.
	\end{lemma}
\noindent 		{\bf Proof:}
	 From Assumption \ref{ass:ass1} we have
\begin{align*}
	\vert f_{yy}(t,y)\vert=& \vert \exp\{-g(y-bt)\}[g'(y-bt)]^2-\exp\{-g(y-bt)\}g''(y-bt)\vert\\
	\leq & (L^2+L)e^{-g(y-bt)}\leq(L^2+L)e^{-g(-bt)} =:c(t)
\end{align*}
for all $y\in[0,\infty)$. Consequently  $y\mapsto f_y(t,y)$ is a Lipschitz function on $\mathbb{R}_+$, uniformly in $t$ in compact intervals. 
Therefore
\begin{eqnarray}\label{new-4-45}
\vert C_n(t)\vert\leq a\int_0^t c(s)\left(n^{\frac{\alpha+\beta-3}{4}}\vert \bar{Y}_n(s)-U(s)\vert\right)^2ds.
\end{eqnarray}
 The result now follows by noting that $\frac{\alpha+\beta-3}{4}<\beta-2$ (see \eqref{eq:ab1210}) and using \eqref{eq:strong117}. $\Box$\\

	\begin{lemma}\label{bnto0}
	Let $D_n$ be as given in \eqref{D-n}. 	As $n\to \infty$, $\sup_{0 \le s \le t} |D_n(s)| \to 0$, in probability, for every $t \ge 0$.
		\end{lemma}
\noindent 		{\bf Proof:}
 Note that
 \begin{eqnarray}\label{s-2-45}
D_n(t)&=&n^{\frac{\alpha-\beta+1}{2}}\left(\int_0^t\frac{a\theta^{2-\beta}f(s,U(s))}{(t-s)^{\beta-2}}ds-n^{\beta-2}\int_0^t\frac{a f(s,\bar{Y}_n(s))}{(n\theta(t-s)+1)^{\beta-2}}ds\right)\nonumber\\
&=&n^{\frac{\alpha-\beta+1}{2}}\left(\int_0^t\frac{a\theta^{2-\beta}f(s,U(s))}{(t-s)^{\beta-2}}ds-n^{\beta-2}\int_0^t\frac{a f(s,U(s))}{(n\theta (t-s)+1)^{\beta-2}}ds\right)\nonumber\\
&+&n^{\frac{\alpha-\beta+1}{2}}\left(n^{\beta-2}\int_0^t\frac{a f(s,U(s))}{(n\theta (t-s)+1)^{\beta-2}}ds-n^{\beta-2}\int_0^t\frac{a f(s,\bar{Y}_n(s))}{(n\theta (t-s)+1)^{\beta-2}}ds\right).\nonumber\\
 \end{eqnarray}
For the first term, note that
 \begin{eqnarray}\label{s-2-46}
 0&<&\int_0^t\frac{a\theta^{2-\beta}f(s,U(s))}{(t-s)^{\beta-2}}ds-n^{\beta-2}\int_0^t\frac{a f(s,U(s))}{(n\theta (t-s)+1)^{\beta-2}}ds\nonumber\\
 &=&\int_0^ta\theta^{2-\beta}f(s,U(s))\left(\frac{1}{(t-s)^{\beta-2}}-\frac{1}{(t-s+\frac{1}{n\theta})^{\beta-2}}\right)ds\nonumber\\
 &\leq&\frac{a\theta^{2-\beta}}{3-\beta}\sup_{0\leq s\leq t}\{\exp\{-g(-bs)\}\}\left( t^{3-\beta}-(t+\frac{1}{n\theta})^{3-\beta}+\left(\frac{1}{n\theta}\right)^{3-\beta}\right)\nonumber\\
 &<&\frac{a\theta^{2-\beta}}{3-\beta}\sup_{0\leq s\leq t}\{\exp\{-g(-bs)\}\}\left(\frac{1}{n\theta}\right)^{3-\beta}.
 \end{eqnarray}
Also, from the Lipschitz property of $f$ (see (\ref{new-3-25})) we have
\begin{eqnarray}\label{s-2-47}
&&\left\vert n^{\beta-2}\int_0^t\frac{a f(s,U(s))}{(n\theta (t-s)+1)^{\beta-2}}ds-n^{\beta-2}\int_0^t\frac{a f(s,\bar{Y}_n(s))}{(n\theta(t-s)+1)^{\beta-2}}ds\right\vert\nonumber\\
&\leq& \int_0^t\frac{a\theta^{2-\beta} Le^{-g(-bs)}\vert \bar{Y}_n(s)-U(s)\vert}{(t-s+\frac{1}{n\theta})^{\beta-2}}ds\nonumber\\
&\leq&\frac{a\theta^{2-\beta}L t^{3-\beta}}{3-\beta}\sup_{0\leq s\leq t}\{\exp\{-g(-bs)\}\}\sup_{0\leq s\leq t}\{\vert \bar{Y}_n(s)-U(s)\vert\}.
\end{eqnarray}
Combining (\ref{s-2-45})-(\ref{s-2-47}), we have
\begin{equation}\label{e-3-45}
\sup_{0\leq s\leq t}\vert D_n(s)\vert\leq C \left [ n^{\frac{\alpha + \beta - 5}{2}}  + n^{\frac{\alpha-\beta+1}{2}}\vert U_n(s)-U(s)\vert\right].
\end{equation}
From \eqref{eq:ab1210} we see that $\frac{\alpha + \beta - 5}{2}< 0$ and $\frac{\alpha-\beta+1}{2}<\beta-2$.  The result follows using \eqref{eq:strong117}. $\Box$\\

\begin{lemma} \label{sicgce}
Let $\mathbf{R}_n=(R_{n,1},\dots,R_{n,d})^T$ and $\mathbf{R}=(R_1,\dots,R_d)^T$ be as given by \eqref{R-n} and \eqref{R} respectively. As $n \to \infty$, $\mathbf{R}_n$ converges to $\mathbf{R}$ in distribution, in $D_{\RR^d}[0,\infty)$.	
\end{lemma}
\noindent 		{\bf Proof:} Let $\tilde{\mathbf{ R}}_n = (\tilde R_{n,1}, \dots, \tilde R_{n,d})^T$, where
$$
\tilde R_{n,i}(t) = \int_{B(t)}[r\wedge(t-\gamma(s))]\Sigma_{n,i}(ds,dr), \; t \ge 0.$$
Since  $\Lambda_n(t)$ is an $\{\mathcal{F}_u^n\}$-stopping time for each $t\geq 0$,  $1_{[0, \Lambda_n(t)]}(s)[r\wedge(t-\gamma_n(s))]$ is $\mathcal{F}_s^n$-predictable. Thus applying the isometry property of the stochastic integral, we obtain 
\begin{eqnarray}\label{e-3-25}
	&&\EE\left(R_{n,i}(t) - \tilde R_{n.i}(t)\right)^2\nonumber\\
&=&\EE\left(\int_{\RR_+\times\RR_+}\left (1_{[0, \Lambda_n(t)]}(s)[r\wedge(t-\gamma_n(s))]
-1_{[0, \Lambda(t)]}(s)[r\wedge(t-\gamma(s))] \right )\Sigma_{n,i}(ds,dr)\right)^2\nonumber\\
&=&\EE\int_0^\infty\int_0^\infty\left(1_{[0, \Lambda_n(t)]}(s)[r\wedge(t-\gamma_n(s))]-1_{[0, \Lambda(t)]}(s)[r\wedge(t-\gamma(s))]\right)^2n^{\beta-1}\nu_n(dr)ds\nonumber\\
&\leq&2\EE\left(\int_0^\infty\int_0^\infty 1_{[0,\Lambda_n(t)]}(s)[r\wedge(t-\gamma_n(s))-r\wedge(t-\gamma(s))]^2n^{\beta-1}\nu_n(dr)ds\right)
  \nonumber\\
&&+2\EE\left(\vert\Lambda_n(t)-\Lambda(t)\vert\right)\,\int_0^\infty (r\wedge t)^2n^{\beta-1}\nu_n(dr)\nonumber\\
&=&2\EE\left(\int_0^\infty\int_0^\infty 1_{[0,\Lambda_n(t)]}(s)[r\wedge(t-\gamma_n(s))-r\wedge(t-\gamma(s))]^2\frac{\theta (\beta-1)n^\beta}{(n\theta r+1)^{\beta}}dr\,ds\right)
  \nonumber\\
&&+2\EE\left(\vert\Lambda_n(t)-\Lambda(t)\vert\right)\,\int_0^\infty (r\wedge t)^2\theta (\beta-1)n^\beta(n\theta r+1)^{-\beta}dr.
\end{eqnarray}                 
By the dominated convergence theorem, we have by using the fact that $\beta \in (2,3)$,
\begin{equation}\label{eq:dct1}
\lim\limits_{n\to\infty}\int_0^\infty (r\wedge t)^2\theta (\beta-1)n^\beta(n\theta r+1)^{-\beta}dr=\int_0^\infty(r\wedge t)^2\theta^{1-\beta} (\beta-1)r^{-\beta}dr<\infty.
\end{equation}

Note that
$$
0 \le \min \{ \Lambda(t), \Lambda_n(t) \} \le \max \{ \Lambda(t), \Lambda_n(t) \} \le \int_0^t \exp\{-g(-bs)\}ds.$$
Consequently
$$
\vert\Lambda_n(t) - \Lambda(t)\vert \le 2 \int_0^t \exp\{-g(-bs)\}ds.$$
Also, from \eqref{new-3-25} we have
$$
\vert\Lambda_n(t) - \Lambda(t)\vert \le L \sup_{0\leq s\leq t}\vert U_n(s)-U(s)\vert \int_0^t \exp\{-g(-bs)\}ds.$$
Thus \eqref{eq:strong117} and the dominated convergence theorem yield
%
%
%
\begin{eqnarray}\label{e-3-27}
\lim_{n\to \infty}\EE\left(\vert\Lambda_n(t)-\Lambda(t)\vert\right) = 0.
\end{eqnarray}
Combining \eqref{eq:dct1} and \eqref{e-3-27} we have that the second term on the right-hand side of \eqref{e-3-25} converges to $0$ as $n\to \infty$.

Now we consider the first term.
From the definitions of $\gamma_n$ and $\gamma$ we see, for any $s \ge 0$,
\[
s=\int_0^{\gamma_n(s)}f(z,\bar{Y}_n(z)dz=\int_0^{\gamma(s)}f(z,U(z))dz.
\]
Consequently,
\begin{equation}\label{s-s-2-42}
\int_0^{\gamma(s)}f(z,U(z))dz-\int_0^{\gamma_n(s)}f(z,U(z))dz=\int_0^{\gamma_n(s)}f(z,\bar{Y}_n(z))dz-\int_0^{\gamma_n(s)}f(z,U(z))dz. 
\end{equation}
Since $f(z,U(z))$ is bounded below away from $0$ (see Remark \ref{R-3-2}), there exists a  $c>0$ such that 
\begin{equation}\label{s-s-2-43}
\left\vert \int_0^{\gamma(s)}f(z,U(z))dz-\int_0^{\gamma_n(s)}f(z,U(z))dz\right\vert\geq c\vert \gamma(s)-\gamma_n(s)\vert.
\end{equation}
On the other hand, from \eqref{new-3-25}, we obtain that for any $s\leq \Lambda_n(t)$ (equivalently, $\gamma_n(s)\leq t$)
\begin{eqnarray}\label{s-s-2-44}
&&\left\vert \int_0^{\gamma_n(s)}f(z,\bar{Y}_n(z))dz-\int_0^{\gamma_n(s)}f(z,U(z))dz\right\vert\nonumber\\
&\leq&L \sup_{0\leq z\leq \gamma_n(s)}\vert \bar{Y}_n(z)-U(z)\vert \int_0^{\gamma_n(s)}\exp\{-g(-bu)\}du\nonumber\\
&\leq& L\sup_{0\leq s\leq t}\vert \bar{Y}_n(s)-U(s)\vert \int_0^{t}\exp\{-g(-bu)\}du.
\end{eqnarray}
Combining (\ref{s-s-2-42})-(\ref{s-s-2-44}) we have that 
\begin{equation}
1_{[0,\Lambda_n(t)]}(s)\vert \gamma_n(s)-\gamma(s)\vert\leq \frac{L}{c} \sup_{0\leq s\leq t}\vert \bar{Y}_n(s)-U(s)\vert \int_0^{t}\exp\{-g(-bu)\}du.
\end{equation}
Using \eqref{eq:strong117} we now obtain
\begin{equation}
1_{[0,\Lambda_n(t)]}(s)\vert \gamma_n(s)-\gamma(s)\vert\to 0
\end{equation}
in probability, as $n\to \infty$. An application of the dominated convergence theorem now shows that
\begin{equation}\label{e-3-31}
\lim\limits_{n\to\infty}\EE\left(\int_0^\infty\int_0^\infty 1_{[0,\Lambda_n(t)]}(s)[r\wedge(t-\gamma_n(s))-r\wedge(t-\gamma(s))]^2\frac{\theta (\beta-1)n^\beta}{(n\theta r+1)^{\beta}}dr\,ds\right) = 0.
\end{equation}
Thus the first term on the right-hand side of \eqref{e-3-25} converges to $0$ as well. 

Combining the above observations we have that for each $i=1,\dots,d$ and $t \ge 0$
\begin{equation}\label{new-R-4-22}
\lim\limits_{n\to\infty}\EE\left(R_{n,i}(t) - \tilde R_{n.i}(t)\right)^2=0.
\end{equation}
Note that 
\[
\tilde R_{n.i}(t)=\int_{[0,\infty)\times[0,\infty)}n^{\frac{\beta-\alpha-1}{2}}1_{B(t)}(s,r)(r\wedge(t-\gamma(s)))\tilde{\xi}_{n,i}(ds,dr).
\]

For each fixed $i=1,\dots,d$, let us show the weak convergence of the finite-dimensional distribution of $\tilde R_{n,i}$. For any $0<t_1<\cdots<t_k<\infty$, denote $f^n(s,r)=(f_1^n(s,r),\dots,f_k^n(s,r))^T$ where $f_j^n(s,r)=n^\frac{\beta-\alpha-1}{2}1_{B(t_j)}(s,r)(r\wedge(t_j-\gamma(s))$, $j=1,\dots,k$. Then $$\tilde{R}_{n,i}(t_j)=\int_{[0,\infty)\times[0,\infty)}f_j^n(s,r)\tilde{\xi}_{n,i}(ds,dr).$$ One can show by the change of variables that
\begin{eqnarray}
&&\lim\limits_{n\to\infty}\int_{[0,\infty)\times[0,\infty)}f_j^n(s,r)f_l^n(s,r)\eta_n(ds,dr)\nonumber\\
&=&\lim\limits_{n\to\infty}\int_{[0,\infty)\times[0,\infty)}f_j^n(s,r)f_l^n(s,r)n^\alpha (\beta-1)n\theta(n\theta r+1)^{-\beta}dr\,ds\nonumber\\
&=&\lim\limits_{n\to\infty}\int_0^{t_j\wedge t_l}\int_0^\infty f(s,U(s))[r\wedge(t_j-s)][r\wedge(t_l-s)]n^\beta (\beta-1)\theta(n\theta r+1)^{-\beta}dr\,ds\nonumber\\
&=&\int_0^{t_j\wedge t_l}\int_0^\infty f(s,U(s))[r\wedge(t_j-s)][r\wedge(t_l-s)] (\beta-1)\theta^{1-\beta}r^{-\beta}dr\,ds\nonumber\\
&=&\EE(R_i(t_j)R_{i}(t_l)).
\end{eqnarray}
Since $\vert f^n\vert \leq n^\frac{\beta-\alpha-1}{2} t_k$ and $\lim\limits_{n\to\infty}n^{\frac{\beta-\alpha-1}{2}}=0$, we deduce that $1_{\{\vert f^n\vert>\varepsilon\}}=0$ for large enough $n$, and hence, for each $\varepsilon>0$ and $j$, if $n$ is large enough, 
$$\int_{[0,\infty)\times[0,\infty)}1_{\{\vert f^n\vert>\varepsilon\}}\vert f_j^n(s,r)\vert^2\eta_n(ds,dr)=0.
$$
From Theorem 6.1 in \cite{K1} it now follows that $(\tilde{R}_{n,i}(t_1),\dots,\tilde{R}_{n,i}(t_k))^T\Rightarrow (R_i(t_1),\dots,R_i(t_k))^T$ as $n$ increases to infinity, for each $i=1,\dots,d$. Since $\tilde{\mathbf{R}}_n$ has independent components, we have that the  finite-dimensional distributions of $\tilde{\mathbf{R}}_n$ converge to those of $\mathbf{R}$. Using \eqref{new-R-4-22}, we then obtain that the  finite-dimensional distributions of $\mathbf{R}_n$ converge to those of $\mathbf{R}$.


It thus suffices to show that $\{\mathbf{R}_n\}$ is tight in $D_{\RR^d}[0, \infty)$, for which, it suffices to prove the  following estimate: for each fixed $T>0$ there exists a constant $C_T>0$ such that for $0\leq h\leq 1$ and $0\leq t\leq T$
\begin{equation}
\EE\left(\left\vert \mathbf{R}_n(t+h)-\mathbf{R}_n(t)\right\vert^4\right)\leq C_Th^2,
\end{equation}
 Recall the definition of $\clu_{n,i}^{(2)}$ in \eqref{eq:clu3}.
Then
$$R_{n,i}(t+h)-R_{n,i}(t)=n^{\frac{\alpha + \beta -3}{2}}\clu_{n,i}^{(2)}(\Lambda_n(t+h))=n^{\frac{\alpha + \beta -3}{2}}\left(A_{n,i}(t+h)-A_{n,i}(t)\right).$$
 From \eqref{s-2-21}, \eqref{s-2-22} and \eqref{s-2-23}, we now have 
\begin{eqnarray}
\EE\left(\vert \mathbf{R}_n(t+h)-\mathbf{R}_n(t)\vert^4\right)&=&n^{\frac{4(\alpha+\beta-3)}{2}}\EE\left(\vert \mathbf{A}_n(t+h)-\mathbf{A}_n(t)\vert^4\right)\nonumber\\
&\leq&Cn^{2\alpha+2\beta-6}\left(n^{5-\beta-3\alpha}h^{6-\beta}+n^{6-2\beta-2\alpha}h^{2(4-\beta)}\right)\nonumber\\
&=&C\left(n^{-(\alpha-\beta+1)}h^{6-\beta}+h^{2(4-\beta)}\right)\leq Ch^2,
\end{eqnarray}
where the last inequality follows from $\alpha>\beta-1$ and $2<\beta<3$.
This proves the desired tightness and the result follows. $\Box$

%
%

\bigskip

\setcounter{equation}{0}
\section{The Moment Stabilization Property }\label{sec:sec5}
In this section, we will prove Theorem \ref{momstab}.
Let $\mathbf{Z}$ be as in Theorem \ref{CLT} and let
 $\mathbf{R}=(R_1,\dots, R_d)^T$ be the Gaussian process introduced in \eqref{R}.
Then  $\bar{Z}=\frac{1}{d}\sum_{i=1}^d Z_i$ satisfies 
\begin{equation}\label{bar}
 \bar{Z}(t)= \bar R(t) + a\int_0^tf_y(s,U(s))\bar{Z}(s)ds
\end{equation}
where
$
\bar R =\frac{1}{d}\sum_{i=1}^dR_i
$.
 Note that $\bar R$ is a zero mean Gaussian process. We begin by computing the covariance functions of $R_i$, $i=1,\dots, d$, and $\bar{R}$.
\begin{lemma}\label{L-5-1} 
The covariance functions of the Gaussian processes $R_i$, $i=1,\dots,d$, and $\bar R$ are given respectively by
\begin{equation}\label{e-5-3-0}
cov( R_i(s), R_i(t))=\mathbb{E}( R_i(s) R_i(t))=\theta^{1-\beta}\int_0^s\int_0^t\int_0^{u\wedge v}\exp\{-g(U(z)-bz)\}(u\vee v-z)^{1-\beta}dz\,du\,dv, 
\end{equation}
and
\begin{equation}\label{e-5-3}
cov(\bar R(s),\bar R(t))=\mathbb{E}(\bar R(s)\bar R(t))=\frac{\theta^{1-\beta}}{d}\int_0^s\int_0^t\int_0^{u\wedge v}\exp\{-g(U(z)-bz)\}(u\vee v-z)^{1-\beta}dz\,du\,dv, 
\end{equation}
for any $s, t\geq0$.
\end{lemma}
\noindent \textbf{Proof:}
Without loss of generality, we assume that  $0\leq s\leq t<\infty$. For each $i=1,\dots,d$, by the change of variables $u=\int_0^z \exp\{-g(U(x)-bx)\}dx$, we have
\begin{eqnarray}\label{new-5-3}
&&\theta^{\beta-1}\EE(R_i(s)R_i(t))\nonumber\\
&=&\int_{B(s)}[r\wedge(s-\gamma(u))][r\wedge(t-\gamma(u))](\beta-1)r^{-\beta}dr\,du\nonumber\\
&=&\int_0^s\int_0^\infty \exp\{-g(U(z)-bz)\}[r\wedge(s-z)][r\wedge(t-z)](\beta-1)r^{-\beta}dr\,dz\nonumber\\
&=&\int_0^s\int_0^{s-z}\exp\{-g(U(z)-bz)\}r^2(\beta-1)r^{-\beta}dr\,dz\nonumber\\
&&+\int_0^s\int_{s-z}^{t-z}\exp\{-g(U(z)-bz)\}(s-z)r(\beta-1)r^{-\beta}dr\,dz\nonumber\\
&&+\int_0^s\int_{t-z}^\infty \exp\{-g(U(z)-bz)\}(s-z)(t-z)(\beta-1)r^{-\beta}dr\,dz\nonumber\\
&=&\int_0^s\exp\{-g(U(z)-bz)\}\nonumber\\
&&\quad \times\left[\frac{2}{3-\beta}(s-z)^{3-\beta}+\frac{1}{\beta-2}(s-z)[(s-z)^{2-\beta}-(t-z)^{2-\beta}]\right]dz.\nonumber\\
\end{eqnarray}
A simple calculation shows that
\begin{equation}\label{new-5-4}
\frac{2}{3-\beta}(s-z)^{3-\beta}+\frac{1}{\beta-2}(s-z)[(s-z)^{2-\beta}-(t-z)^{2-\beta}]=\int_z^s\int_z^t (u\vee v-z)^{1-\beta}du\,dv.
\end{equation}
Combining (\ref{new-5-3}) and (\ref{new-5-4}) and by changing the order of integration, we obtain 
\begin{eqnarray*}\label{covariance-1}
\theta^{\beta-1}\EE(R_i(s)R_i(t))
&=&\int_0^s\exp\{-g(U(z)-bz)\}\int_z^s\int_z^t (u\vee v-z)^{1-\beta}du\,dv\,dz\nonumber\\
&=&\int_0^s\int_0^t\int_0^{u\wedge v}\exp\{-g(U(z)-bz)\}(u\vee v-z)^{1-\beta}dz\,du\,dv.
\end{eqnarray*}
This proves \eqref{e-5-3-0}.
Equation \eqref{e-5-3} is now immediate on noting that $R_1, \dots, R_d$ are i.i.d.\ $\Box$

\bigskip

In the next lemma, we give a bound on the second moments of the increment of the Gaussian processes $R_i$, $i=1,\dots, d$, and $\bar{R}$.
\begin{lemma}\label{lem-5-2}
For any $s,t\geq 0$, the following bounds hold:
\begin{equation}\label{R-bound}
\EE(\vert  R_i(t)-R_i(s)\vert^2)\leq \frac{2K_1\theta^{1-\beta}}{(\beta-2)(3-\beta)(4-\beta)}(t-s)^{4-\beta},
\end{equation}
and
\begin{equation}\label{bar-R-bound}
\EE(\vert \bar R(t)-\bar R(s)\vert^2)\leq \frac{2K_1\theta^{1-\beta}}{d(\beta-2)(3-\beta)(4-\beta)}(t-s)^{4-\beta},
\end{equation}
where $K_1:=\sup_{s\geq 0}\{\exp\{-g(U(s)-bs)\}\}$. 

Consequently, Gaussian processes $R_1,\dots, R_d, \bar{R}$ have  versions that are H\"{o}lder continuous  of any order $\rho\in (0,(4-\beta)/2)$, on $[0,T]$, for all $T > 0$. 
\end{lemma}
\noindent \textbf{Proof:}
Fix $0 \le s\leq t < \infty$. From  Lemma \ref{L-5-1}, for each $i=1,\dots,d$,
\begin{eqnarray}
\EE(\vert  R_i(t)-R_i(s)\vert^2)&=&cov(  R_i(t), R_i(t))-2cov(  R_i(s), R_i(t))+cov(  R_i(s), R_i(s))\nonumber\\
&=&\theta^{1-\beta}\int_s^t\int_s^t\int_0^{u\wedge v}\exp\{-g(U(z)-bz)\}(u\vee v-z)^{1-\beta}dz\,du\,dv\nonumber\\
&\leq&K_1\theta^{1-\beta}\int_s^t\int_s^t\int_0^{u\wedge v}(u\vee v-z)^{1-\beta}dz\,du\,dv\nonumber\\
&=&\frac{K_1\theta^{1-\beta}}{\beta-2}\int_s^t\int_s^t\left[ \vert u-v\vert^{2-\beta}-(u\vee v)^{2-\beta}\right]du\,dv\nonumber\\
&<&\frac{K_1\theta^{1-\beta}}{\beta-2}\int_s^t\int_s^t \vert u-v\vert^{2-\beta}du\,dv\nonumber\\
&=&\frac{2K_1\theta^{1-\beta}}{\beta-2}\int_s^t\int_v^t  (u-v)^{2-\beta}du\,dv\nonumber\\
&=&\frac{2K_1\theta^{1-\beta}}{(\beta-2)(3-\beta)(4-\beta)}(t-s)^{4-\beta}.\label{eq:5.9}
\end{eqnarray}
This completes the proof of \eqref{R-bound}. Inequality in \eqref{bar-R-bound} is now immediate.  The second statement in the lemma now follows from Kolmogorov's continuity  criterion. $\Box$\\ \\

Proof of Theorem \ref{momstab} relies on an explicit representation for the solution of equation \eqref{bar}. For that we begin with an indefinite integral of a deterministic function with respect to the Gaussian process $\bar{R}$. 

Denote by $\mathcal{E}$ the linear span of indicator functions of the form $ 1_{(s,t]}: \RR_+ \to \RR$, $0 \le s \le t < \infty$. Consider the inner product on $\mathcal{E}$ given by 
\begin{equation}
\langle 1_{(0,s]}, 1_{(0,t]}\rangle_{\mathcal{H}_{\bar{R}}}=cov(\bar R(s),\bar R(t)).
\end{equation}
We denote by $\mathcal{H}_{\bar{R}}$ the Hilbert space obtained as the closure of $\mathcal{E}$ with respect to this inner product. 
Define $\bar\bfr: \cle \to L^2(\Omega, \clf, \PP)$ as
$$
\bar\bfr(1_{(0,t])}) = \bar R(t),\ 0 \le t < \infty, 
$$
where the definition is extended to all of $\cle$ by linearity.  Clearly $\EE(|\bar \bfr(\phi)|^2) = \langle \phi, \phi\rangle_{\mathcal{H}_{\bar{R}}}$ 
for all $\phi \in \cle$.  We can now extend the definition of $\bar \bfr$ to all of $\mathcal{H}_{\bar{R}}$ by isometry. Occasionally, we will use the notation
$$\bar \bfr(\phi) = \int_0^{\infty} \phi(t) d\bar R(t), \; \phi \in \mathcal{H}_{\bar{R}}.$$
Define 
\begin{equation}
\rho(u,v)=
\begin{cases}
\frac{\theta^{1-\beta}}{d}\int_0^{u\wedge v}\exp\{-g(U(z)-bz)\}(u\vee v-z)^{1-\beta}dz, \,\hbox{\text{if}}\, u\neq v;\\
0,\,\hbox{\text{if}}\  u=v,
\end{cases}
\end{equation}
 for all $u, v\geq 0$. From Lemma \ref{L-5-1}, we see that
\begin{equation}\label{new-5-7}
cov(\bar{R}(s),\bar{R}(t))=\int_0^s\int_0^t\rho(u,v)du\,dv
\end{equation}
and for any $\phi,\ \tilde \phi \in \cle$
\begin{equation}\label{eq:phiisom}
\langle \phi, \tilde \phi\rangle_{\mathcal{H}_{\bar{R}}}=\int_0^\infty \int_0^\infty \phi(u)\tilde \phi(v)\rho(u,v)du\,dv.
\end{equation}
It can be shown that  $\mathcal{H}_{\bar{R}}$ contains all measurable functions $\phi$ on $\RR_+$ satisfying 
\begin{equation}\label{e-5-11}
\int_0^\infty\int_0^\infty \vert\phi(u)\vert \vert\phi(v)\vert \rho(u,v)du\,dv<\infty.
\end{equation}
One can also show that equality \eqref{eq:phiisom} holds for $\phi,\ \tilde{\phi}$ that satisfy \eqref{e-5-11}.

This type of isometry is considered in \cite{N06} (see Chapter 5) and \cite{AMN} with respect to fractional Brownian motion and general Gaussian processes respectively. 



\begin{remark}
If $\phi: [0,\infty)\to\mathbb{R}$ is continuous, then, for any $t>0$, the function $\phi_t$ defined by $\phi_t(\cdot)=1_{[0,t]}(\cdot)\phi(\cdot)$ satisfies (\ref{e-5-11}).  Consequently,
 $\phi_t$ is in $\mathcal{H}_{\bar{R}}$ and we write formally
\begin{equation}\label{new-5-17}
\bar \bfr(\phi)(t):= \bar \bfr(\phi_t) = \int_0^t\phi(s)d\bar R(s).
\end{equation}
\end{remark}

\begin{remark}\label{R-5-4}
If $\phi(\cdot)$ is H\"{o}lder continuous of order $\rho_1>1-\frac{4-\beta}{2}$ on $[0,t]$, for every $t>0$, as a result of Young integration theory \cite{Y}, the pathwise Riemann-Stieltjes integral $\int_0^t\phi(s)d\bar{R}(s)$ exists, since $\bar{R}$ is H\"{o}lder continuous of any order $\rho\in (0,\frac{4-\beta}{2})$.  
Z\"{a}hle \cite{Z}  showed
(see Proposition 4.4.1 therein) that $\bar \bfr(\phi)(t)$ is H\"{o}lder continuous of the same order as $\bar{R}$ on $[0,T]$, for every $T>0$. The indefinite integral $\bar \bfr(\phi)(\cdot)$ on the right-hand side of \eqref{new-5-17} coincides with the pathwise Riemann-Stieltjes integral.
%
\end{remark}

We now proceed to the proof of Theorem \ref{momstab}.\\

\noindent \textbf{Proof of Theorem \ref{momstab}:}
 Define
 \[
\phi(t):=\exp\left\{-a\int_0^tf_y(z,U(z))dz\right\},
\;
   \tilde{\phi}(t)=\exp\left\{a\int_0^t f_y(z,U(z))dz\right\}.\]
Then the derivatives of $\phi$ and $\tilde \phi$ are 
\[
\phi'(t)=-af_y(t,U(t))\phi(t),
\;
 \tilde{\phi}'(t)=af_y(t,U(t)) \tilde{\phi}(t).
\]
Remark \ref{R-3-2} implies that $\phi'$ and $\tilde{\phi}'$ are bounded on any compact interval, and hence, $\phi$ and $\tilde{\phi}$ are locally Lipschitz continuous. From Remark \ref{R-5-4}, the indefinite integral 
$$\bar \bfr(\phi)(t)=\int_0^t\phi(s)d\bar{R}(s)=\int_0^t \exp\left\{-a\int_0^s  f_y(z,U(z))dz\right\}d\bar R(s)$$ is well defined as a Riemann-Stieltjes integral, and for every $T>0$,  $\bar \bfr(\phi)(t)$ is H\"{o}lder continuous on $[0,T]$ of any order $\rho \in(0,(4-\beta)/2)$.

It follows from Theorems 3.1 and 4.4.2 in \cite{Z} that 
\begin{eqnarray}\label{inte-by-par}
 \tilde{\phi}(t)\bar \bfr(\phi)(t)&=&\int_0^t \tilde{\phi}(s)\phi(s)d\bar{R}(s)+a\int_0^tf_y(s,U(s))  \tilde{\phi}(s)\bar \bfr(\phi)(s)ds\nonumber\\
&=&\bar{R}(t)+a\int_0^tf_y(s,U(s)) \tilde{\phi}(s)\bar \bfr(\phi)(s)ds,
\end{eqnarray}
which implies that $ \tilde{\phi}(t)\bar \bfr(\phi)(t)$ solves \eqref{bar}. Thus, the solution $\bar{Z}$ to \eqref{bar} can be written explicitly as
\begin{eqnarray}\label{e-5-12}
\bar{Z}(t)&=& \tilde{\phi}(t)\bar \bfr(\phi)(t)=\exp\left\{a\int_0^t f_y(z,U(z))dz\right\}\int_0^t \exp\left\{-a\int_0^s  f_y(z,U(z))dz\right\}d\bar R(s).\nonumber\\
\end{eqnarray}
 By the isometry of the mapping $\bar\bfr$ we have, on letting $\phi_t(u) = \phi(u)1_{[0,t]}(u)$,
\begin{eqnarray}\label{e-5-14}
\EE(\vert\bar{ Z}(t)\vert^2)&=&\vert \tilde{\phi}(t)\vert^2\int_0^\infty\int_0^\infty \phi_t(u)\phi_t(v)\rho(u,v)du\,dv\nonumber\\
&=&\vert \tilde{\phi}(t)\vert^2\int_0^t\int_0^t\phi(u)\phi(v)\rho(u,v)du\,dv\nonumber\\
&=&\int_0^t\int_0^t\exp\left\{a\int_u^tf_y(z,U(z))dz\right\}\exp\left\{a\int_v^tf_y(z,U(z))dz\right\}\rho(u,v)du\,dv.\nonumber\\
\end{eqnarray}
Recall the definition of $\mu$ from the statement of Theorem \ref{momstab}.  From Remark \ref{R-3-2}, $\mu \in (0, \infty)$.
Then
by a calculation similar to \eqref{eq:5.9}, we have
\begin{eqnarray}\label{e-5-15}
\theta^{\beta-1}\EE(\vert \bar{Z}(t)\vert^2)&\leq &\int_0^t\int_0^t e^{-a\mu(t-u)}e^{-a\mu(t-v)}\rho(u,v)du\,dv\nonumber\\
&=&\frac{1}{d}\int_0^t\int_0^t\int_0^{u\wedge v}e^{-a\mu(t-u)}e^{-a\mu(t-v)}(u\vee v-z)^{1-\beta}dz\,du\,dv\nonumber\\
&\le&\frac{2}{d(\beta-2)}\int_0^t\int_v^t e^{-a\mu(t-u)}e^{-a\mu(t-v)} (u-v)^{2-\beta}du\,dv \nonumber\\
&\leq&\frac{2}{d(\beta-2)}\int_0^t\int_v^t e^{-a\mu(t-v)}(u-v)^{2-\beta}du\,dv\nonumber\\
&=&\frac{2}{d(\beta-2)(3-\beta)}\int_{0}^te^{-a\mu(t-v)}(t-v)^{3-\beta}dv\nonumber\\
&\le&\frac{2}{d(\beta-2)(3-\beta)(a\mu)^{4-\beta}}\Gamma(4-\beta),
\end{eqnarray}
for all $t\geq 0$. $\Box$

\bigskip

\setcounter{equation}{0}
\section{Fractional Ornstein-Uhlenbeck Process}

We now proceed to the proof of Theorem  \ref{ouplim} in this section.
 Throughout this section we take $b=a$.
From Proposition \ref{L-3-1} it follows that
\begin{equation}
	\label{eq:beqa}
	U(t) = bt, \; f(t,U(t)) = 1, \; f_y(t, U(t)) = - g'(0), \mbox{ for all } t \ge 0.
	\end{equation}
For notational simplicity, we will only present the proof for the case $\theta = 1$.
 %

In this special case, the SDE \eqref{bar} can be written as
\begin{equation}\label{e-3-53}
\bar{Z}(t)=\bar{R}(t)-a\int_0^tg'(0)\bar{Z}(s)ds = \bar{R}(t)-\kappa\int_0^t\bar{Z}(s)ds,
\end{equation}
where $\kappa = a\mu = a g'(0)>0$, and 
\[
\bar{R}(t)=\frac{1}{d}\sum_{i=1}^d\int_0^t\int_0^\infty r\wedge(t-s)\Sigma_i(ds,dr),
\]
 for any $t\geq 0$.   

From Lemma \ref{L-5-1}, we have that the covariance of $\bar{R}$  is given by
 \begin{eqnarray}\label{e-6-3}
 cov(\bar{R}(s),\bar{R}(t))=\EE(\bar{R}(s)\bar{R}(t))&=&\frac{1}{d}\int_0^s\int_0^t\int_0^{u\wedge v}(u\vee v-z)^{1-\beta}dz\,du\,dv\nonumber\\
 &=&\frac{1}{d(\beta-2)}\int_0^s\int_0^t(\vert u-v\vert^{2-\beta}-(u\vee v)^{2-\beta})du\,dv,\nonumber\\ 
 \end{eqnarray}
and from Lemma \ref{lem-5-2}, we recall
that the sample paths of the process $\bar R$ are H\"{o}lder continuous on $[0,T]$ of order $\rho$, for any $\rho\in(0,\frac{4-\beta}{2})$.
 
 Recalling  the definition of the indefinite integrals with respect to the Gaussian process $\bar{R}$,  the solution of the SDE in   (\ref{e-3-53}) can be  explicitly written as 
\begin{equation}\label{s-s-3-4}
\bar{Z}(t)=e^{-\kappa t} \int_0^t e^{\kappa s} d\bar{R}(s).
\end{equation}
We now consider the asymptotic behavior of the process $\bar{Z}(t)$, as $t\to \infty$. For  $T,\ t\ge  0$, let $\bar{R}_T(t)=\bar{R(}T+t)-\bar{R}(T)$. 
From (\ref{e-3-53}), we can write 
\begin{equation}\label{e-4-4}
\bar{Z}(T+t)=\bar{Z}(T)+\bar{R}_T(t)-\kappa\int_0^t \bar{Z}(T+s)ds.
\end{equation}
Recall the parameters $\sigma_0^2$ and $\sigma$ introduced before Theorem \ref{ouplim} (also recall that $\theta=1$).
\begin{lemma}\label{L-4-1}
We have
\begin{equation} 
\lim\limits_{T\to\infty}\EE(\vert \bar{Z}(T)\vert^2)=\sigma_0^2.
\end{equation}
\end{lemma}
\noindent \textbf{Proof:}
From (\ref{e-6-3}), (\ref{s-s-3-4}) and the isometry property of of the mapping $\bar\bfr$, we have
\begin{eqnarray}\label{e-6-7}
\EE(\vert \bar{Z}(T)\vert^2)&=&\frac{1}{d\,(\beta-2)}e^{-2\kappa T}\int_{0}^T\int_0^T e^{\kappa v}e^{\kappa u}\left[\vert u-v\vert^{2-\beta}-(u\vee v)^{2-\beta}\right]du\,dv.
\end{eqnarray}
From  (\ref{e-5-15}), 
\begin{equation}
h_1(T):=\frac{1}{d\,(\beta-2)}e^{-2\kappa T}\int_{0}^T\int_0^T e^{\kappa v}e^{\kappa u}\vert u-v\vert^{2-\beta}dudv\leq\frac{2\Gamma(4-\beta)}{d\,(\beta-2)(3-\beta)\kappa^{4-\beta}},
\end{equation}
for all $T\geq 0$.
Also, by  a change of variables, we see that
\begin{eqnarray}
h_1(T)
&=&\frac{1}{d\,(\beta-2)}\int_{0}^T\int_0^Te^{-\kappa v}e^{-\kappa u}\vert u-v\vert^{2-\beta}du\,dv.
\end{eqnarray}
Since the integrand is non-negative, the function $h_1(T)$ is increasing in $T$. Together with the boundedness of $h_1$, the limit $\lim\limits_{T\to\infty}h_1(T)$ exists and is finite.  Furthermore,
\begin{equation}\label{eq:sigz2}
 \lim\limits_{T\to\infty}h_1(T)=\frac{1}{d\,(\beta-2)}\int_{0}^\infty\int_0^\infty e^{-\kappa v}e^{-\kappa u}\vert u-v\vert^{2-\beta}du\,dv = \sigma_0^2.
\end{equation}
Next, for the second term on the right-hand side of (\ref{e-6-7}), note that
\begin{eqnarray}
h_2(T)&:=&\frac{1}{d\,(\beta-2)}e^{-2\kappa T}\int_{0}^T\int_0^T e^{\kappa v}e^{\kappa u}(u\vee v)^{2-\beta}du\,dv\nonumber\\
&=&\frac{2}{d\,(\beta-2)}e^{-2\kappa T}\int_0^T\int_0^ve^{\kappa v}e^{\kappa u}v^{2-\beta}du\,dv\nonumber\\
&=&\frac{2}{d\,\kappa(\beta-2)}e^{-2\kappa T}\int_0^T[e^{2\kappa v}-e^{\kappa v}]v^{2-\beta}du\,dv.
\end{eqnarray}
An application of l'H\^{o}pital's rule shows that
\begin{equation}\label{new-6-12}
\lim\limits_{T\to\infty} h_2(T)=\frac{2}{d\,\kappa(\beta-2)}\lim\limits_{T\to\infty} \frac{ [e^{2\kappa T}-e^{\kappa T}]T^{2-\beta}}{2\kappa e^{2\kappa T}}=0.
\end{equation}
Using \eqref{eq:sigz2} and \eqref{new-6-12}, we now have that
\begin{equation}
\lim\limits_{T\to\infty}\EE(\vert \bar{Z}(T)\vert^2)=\frac{1}{d\,(\beta-2)}\int_{0}^\infty\int_0^\infty e^{-\kappa v}e^{-\kappa u}\vert u-v\vert^{2-\beta}du\,dv =\sigma_0^2.
\end{equation}
The result follows. $\Box$\\

In the next lemma, we compute the limits of second moments and the covariances of $\bar{R}_T(\cdot)$, as $T\to \infty$.
\begin{lemma} \label{L-4-2} 
For any $t \ge 0$,
\begin{equation}\label{eq:eq1012}
\lim\limits_{T\to\infty}\EE(\vert \bar R_T(t)\vert^2)=\frac{2t^{4-\beta}}{d\,(\beta-2)(3-\beta)(4-\beta)}=\sigma^2 t^{4-\beta},
\end{equation}
and, for any $t\geq s\geq 0$,
\begin{equation}\label{eq:eq1015}
\lim\limits_{T\to\infty}\EE(\bar R_T(t)\bar R_T(s))=\frac{\sigma^2}{2}\left[t^{4-\beta}+s^{4-\beta}-(t-s)^{4-\beta}\right].
\end{equation}
\end{lemma}
\noindent \textbf{Proof:} For any $t,\, T\geq 0$, note that 
\[
\bar R_T(t)=\bar R(T+t)-\bar R(t)=\int_0^{T+t}1_{(T,T+t]}(s)d\bar R(s).
\]
By the isometry property of the mapping $\bar\bfr$ and a change of variables, we have 
\begin{eqnarray*}
&&\lim\limits_{T\to\infty}\EE(\vert  \bar R_T(t)\vert^2)\nonumber\\
&=&\lim\limits_{T\to\infty}\frac{1}{d\,(\beta-2)}\int_T^{T+t}\int_T^{T+t}\left(|u-v|^{2-\beta}-(u\vee v)^{2-\beta}\right)du\,dv\nonumber\\
&=&\lim\limits_{T\to\infty}\frac{2}{d\,(\beta-2)}\int_T^{T+t}\int_T^{v}\left((v-u)^{2-\beta}- v^{2-\beta}\right)du\,dv\nonumber\\
&=&\lim\limits_{T\to\infty}\frac{2}{d\,(\beta-2)(3-\beta)}\int_T^{T+t}\left((v-T)^{3-\beta}-(3-\beta) v^{2-\beta}t\right)dv\nonumber\\
&=&\lim\limits_{T\to\infty}\frac{2}{d\,(\beta-2)(3-\beta)}\int_0^{t}\left(v^{3-\beta}- (3-\beta)(v+T)^{2-\beta}t\right)dv\nonumber\\
&=&\frac{2}{d\,(\beta-2)(3-\beta)}\int_0^{t}v^{3-\beta}dv\nonumber\\
&=&\frac{2t^{4-\beta}}{d\,(\beta-2)(3-\beta)(4-\beta)}=\sigma^2t^{4-\beta}.
\end{eqnarray*}
This proves \eqref{eq:eq1012}.

For any $t\geq s\geq 0$, by the isometry property of the mapping $\bar\bfr$ and a change of variables, we have
\begin{eqnarray*}
&&\lim\limits_{T\to\infty}\EE(\bar R_T(t)\bar R_T(s))=\lim\limits_{T\to\infty}\EE((\bar R(T+t)-\bar R(T))(\bar R(T+s)-\bar R(T))\nonumber\\
&=&\lim\limits_{T\to\infty}\frac{1}{d\,(\beta-2)}\int_T^{T+t}\int_T^{T+s}\left(|u-v|^{2-\beta}-(u\vee v)^{2-\beta}\right)du\,dv\nonumber\\
&=&\lim\limits_{T\to\infty}\frac{1}{d\,(\beta-2)}\int_T^{T+s}\int_T^{T+s}\left(|u-v|^{2-\beta}-(u\vee v)^{2-\beta}\right)du\,dv\nonumber\\
&&+\lim\limits_{T\to\infty}\frac{1}{d\,(\beta-2)}\int_{T+s}^{T+t}\int_T^{T+s}\left((v-u)^{2-\beta}-v^{2-\beta}\right)du\,dv\nonumber\\
&=&\lim\limits_{T\to\infty}\EE(\vert  \bar R_T(s)\vert^2)\nonumber\\
&&+\lim\limits_{T\to\infty}\frac{1}{d\,(\beta-2)(3-\beta)}\int_{T+s}^{T+t}\left((v-T)^{3-\beta}-(v-T-s)^{3-\beta}-(3-\beta)v^{2-\beta}s\right)dv\nonumber\\
&=&\frac{2s^{4-\beta}}{d\,(\beta-2)(3-\beta)(4-\beta)}+\frac{1}{d\,(\beta-2)(3-\beta)}\int_{s}^{t}\left(v^{3-\beta}-(v-s)^{3-\beta}\right)dv\nonumber\\
&=&\frac{1}{d\,(\beta-2)(3-\beta)(4-\beta)}\left[s^{4-\beta}+t^{4-\beta}-(t-s)^{4-\beta}\right]\nonumber\\
&=&\frac{\sigma^2}{2}\left[s^{4-\beta}+t^{4-\beta}-(t-s)^{4-\beta}\right],
\end{eqnarray*}
where the fourth equality follows from \eqref{eq:eq1012} and a change of variables. This proves \eqref{eq:eq1015} and the result follows.
$\Box$\\

The next lemma gives the limit of the covariance function of $\bar Z(T)$ and $\bar R_T(t)$, as $T\to\infty$.
\begin{lemma} \label{L-6-3}
For any $t\geq 0$,
\begin{equation}
\lim\limits_{T\to\infty}cov(\bar Z(T),\bar R_T(t))=\lim\limits_{T\to\infty}\EE(\bar Z(T)\bar R_T(t))=\frac{1}{d\,(\beta-2)}\int_0^t\int_0^\infty e^{-\kappa v}(u+v)^{2-\beta}dv\,du.
\end{equation}
\end{lemma}
\noindent \textbf{Proof:} For any $t, \ T\geq 0$, note that $$\bar R_T(t)=\bar R(T+t)-\bar R(T)=\int_0^{T+t}1_{(T,T+t]}(s)d\bar R(s).$$ Then, it follows from the isometry property of the mapping $\bar\bfr$ and \eqref{s-s-3-4} that
\begin{eqnarray}\label{lim-cov-1}
&&\lim\limits_{T\to\infty}cov(\bar Z(T),\bar R_T(t))=\lim\limits_{T\to\infty}\EE(\bar Z(T)\bar R_T(t))\nonumber\\
&=&\lim\limits_{T\to\infty}\frac{e^{-\kappa T}}{d\,(\beta-2)}\int_0^T\int_{T}^{T+t}e^{\kappa v}\left ((u-v)^{2-\beta}-u^{2-\beta}\right)du\,dv\nonumber\\
&=&\lim\limits_{T\to\infty}\frac{e^{-\kappa T}}{d\,(\beta-2)}\int_0^T\int_{T}^{T+t}e^{\kappa v}(u-v)^{2-\beta}du\,dv\nonumber\\
&&-\lim\limits_{T\to\infty}\frac{e^{-\kappa T}}{d\,(\beta-2)}\int_0^T\int_{T}^{T+t}e^{\kappa v}u^{2-\beta}du\,dv.
\end{eqnarray}
For the second term on the right-hand side of \eqref{lim-cov-1}, we have
\begin{eqnarray}\label{lim-cov-2}
\lim\limits_{T\to\infty}\frac{e^{-\kappa T}}{d\,(\beta-2)}\int_0^T\int_{T}^{T+t}e^{\kappa v}u^{2-\beta}du\, dv&\leq&\lim\limits_{T\to\infty}\frac{te^{-\kappa T}}{d\,(\beta-2)}\int_0^T e^{\kappa v}T^{2-\beta}dv\nonumber\\
&=&\lim\limits_{T\to\infty}\frac{t(1-e^{-\kappa T})}{d\,\kappa(\beta-2)T^{\beta-2}}=0.
\end{eqnarray}
For the first term on the right-hand side of \eqref{lim-cov-1}, using a change of variables, we obtain
\begin{eqnarray}\label{lim-cov-3}
&&\lim\limits_{T\to\infty}\frac{e^{-\kappa T}}{d\,(\beta-2)}\int_0^T\int_{T}^{T+t}e^{\kappa v}(u-v)^{2-\beta}du\, dv\nonumber\\
&=&\lim\limits_{T\to\infty}\frac{e^{-\kappa T}}{d\,(\beta-2)(3-\beta)}\int_0^Te^{\kappa v}[(T+t-v)^{3-\beta}-(T-v)^{3-\beta}]dv\nonumber\\
&=&\lim\limits_{T\to\infty}\frac{1}{d\,(\beta-2)(3-\beta)}\int_0^Te^{-\kappa v}[(t+v)^{3-\beta}-v^{3-\beta}]dv\nonumber\\
&=&\frac{1}{d\,(\beta-2)(3-\beta)}\int_0^\infty e^{-\kappa v}[(t+v)^{3-\beta}-v^{3-\beta}]dv\nonumber\\
&=&\frac{1}{d\,(\beta-2)}\int_0^\infty\int_0^t e^{-\kappa v}(u+v)^{2-\beta}du\,dv\nonumber\\
&=&\frac{1}{d\,(\beta-2)}\int_0^t\int_0^\infty e^{-\kappa v}(u+v)^{2-\beta}dv\,du.
\end{eqnarray}
The result follows on combining  \eqref{lim-cov-1}-\eqref{lim-cov-3}. $\Box$\\

\noindent \textbf{Proof of Theorem \ref{ouplim}:}  Define $\varphi: C_{\RR}[0,\infty)\to C_{\RR}[0,\infty)$ by
\[
[\varphi(x)](t)=x(t)-\kappa\int_0^t[\varphi(x)](s)ds.
\]
Then, $\varphi$ is a continuous mapping from $C_{\RR}[0,\infty)$ to $C_{\RR}[0,\infty)$.

For any $t,\, T\geq 0$, denote $\bar{\mathcal{R}}_T(t)=\bar Z(T)+\bar R_T(t)$. Then $\bar{\mathcal{R}}_T$ is a Gaussian process with continuous trajectories, and $\bar Z(T+\cdot)=\varphi(\bar{\mathcal{R}}_T)(\cdot)$. Therefore, in order to prove
that $\bar Z(T+ \cdot)$  converges in distribution, in $C_{\RR}[0,\infty)$ to $Z_{\infty}$, it suffices to show  the  convergence of $\bar{\mathcal{R}}_T$  to $Z_\infty(0)+\sigma B_{H}(\cdot)$, where $Z_\infty(0), \, B_H $ and $\sigma$ are as defined before Theorem \ref{ouplim}. 
From Lemmas \ref{L-4-1}-\ref{L-6-3}, it follows that the finite dimensional distributions of $\bar{\mathcal{R}}_T$ converge to those of $Z_\infty(0)+\sigma B_{H}(\cdot)$.  It thus suffices to verify that
	$\{\bar{\mathcal{ R}}_T(\cdot)\}_{T > 0}$ is tight in $C_{\RR}[0, \infty)$.
	By the Cauchy-Schwartz inequality and \eqref{bar-R-bound}, it follows that  for any $h\geq 0$, $t\geq h$ and $T \ge 0$,
\begin{eqnarray}\label{e-6-12}
&&\EE(\vert \bar{\mathcal{ R}}_T(t+h)-\bar{\mathcal{ R}}_T(t)\vert \vert\bar{ \mathcal{ R}}_T(t)-\bar{\mathcal{ R}}_T(t-h)\vert)\nonumber\\
&=&\EE(\vert \bar R(T+t+h)-\bar R(T+t)\vert \vert \bar R(T+t)-\bar R(T+t-h)\vert)\nonumber\\
&\le& (\EE\vert \bar R(T+t+h)-\bar R(T+t)\vert^2)^{1/2} (\EE\vert \bar R(T+t)-\bar R(T+t-h)\vert^2)^{1/2}\nonumber\\
&\leq& \frac{2h^{4-\beta}}{d\,(\beta-2)(3-\beta)(4-\beta)}.
\end{eqnarray}
Note that $4-\beta > 1$, since $\beta<3$.  The desired tightness now follows from standard results (cf. Theorems 3.8.6 and 3.8.8 in  \cite{EK}). 
This proves the convergence of $\bar Z(T+ \cdot)$ to $Z_{\infty}(\cdot)$.  Stationarity of $Z_{\infty}$ is now immediate.  The result follows. 
$\Box$

\bigskip

\setcounter{equation}{0}
\appendix
\numberwithin{equation}{section}
\section{Auxiliary Results}


Recall that $2<\beta<3$ and $\alpha>\beta-1$.  Also recall the notation $\vartheta_n^{h,t}(r,s)$ and $\tilde\vartheta_n^{h,t}(r,s)$ introduced in the proof of 
Theorem \ref{conv-prob}.  The following lemma gives a key estimate for the proof of the theorem.  
\begin{lemma}\label{L-A-1}
There exists a constant $C>0$ depending only on $\beta$ and $\theta$, such that for any $0\leq h\leq 1$ and $h\leq t<\infty$, we have the following estimates
\begin{equation}\label{s-4-1}
\frac{1}{n^{4(\alpha-1)}} \int_{0}^{t+h}\int_0^\infty  [\tilde \vartheta_n^{h,t}(r,s)]^4n^{\alpha+1} (n\theta r+1)^{-\beta}dr\,ds\leq C n^{5-\beta-3\alpha}h^{6-\beta},
\end{equation}
and
\begin{equation}\label{s-4-2}
\frac{1}{n^{4(\alpha-1)}}\left(\int_{0}^{t+h}\int_0^\infty [\tilde \vartheta_n^{h,t}(r,s)]^2n^{\alpha+1}(n\theta r+1)^{-\beta}dr\,ds\right)^{2}\leq Cn^{6-2\beta-2\alpha}h^{2(4-\beta)}.
\end{equation}
\end{lemma}
\noindent{\textbf{Proof:} We first prove (\ref{s-4-1}). We can write the left-hand side of \eqref{s-4-1} as
\begin{eqnarray}\label{s-4-3}
&&\frac{1}{n^{4(\alpha-1)}}\int_t^{t+h}\int_0^\infty \left[r\wedge(t+h-s)\right]^4n^{\alpha+1}(n\theta r+1)^{-\beta}dr\,ds\nonumber\\
&+&\frac{1}{n^{4(\alpha-1)}}\int_0^{t}\int_0^\infty \left[ r\wedge(t+h-s)-r\wedge(t-s)\right]^4n^{\alpha+1}(n\theta r+1)^{-\beta}dr\,ds.
\end{eqnarray}
We now bound the two terms in (\ref{s-4-3}) separately. In the following, we use $C>0$ to denote a generic constant depending only on $\beta$ and $\theta$; the
value of $C$ may change from one line to next.
Using $(n\theta r+1)^{-\beta}<(n\theta r)^{-\beta}$, we have
\begin{eqnarray}\label{s-4-4}
&&\frac{1}{n^{4(\alpha-1)}}\int_t^{t+h}\int_0^\infty \left[ r\wedge(t+h-s)\right]^4n^{\alpha+1}(n\theta r+1)^{-\beta}dr\,ds\nonumber\\
&=&\frac{1}{n^{4(\alpha-1)}}\int_t^{t+h}\int_0^{t+h-s} r^4n^{\alpha+1}(n\theta r+1)^{-\beta}dr\,ds\nonumber\\
&&+\frac{1}{n^{4(\alpha-1)}}\int_t^{t+h}\int_{t+h-s}^\infty (t+h-s)^4n^{\alpha+1}(n\theta r+1)^{-\beta}dr\,ds\nonumber\\
&\leq&\frac{n^{\alpha-\beta+1}}{n^{4(\alpha-1)}}\int_t^{t+h}\int_0^{t+h-s}r^{4-\beta}drds+\frac{n^{\alpha-\beta+1}}{n^{4(\alpha-1)}}\int_t^{t+h}\int_{t+h-s}^\infty (t+h-s)^4 r^{-\beta}dr\,ds\nonumber\\
&=&C n^{5-\beta-3\alpha}\int_t^{t+h} (t+h-s)^{5-\beta}ds\le Cn^{5-\beta-3\alpha}h^{6-\beta}.
\end{eqnarray} 
For the second term in (\ref{s-4-3}), we have
\begin{eqnarray}\label{s-4-5}
&&\frac{1}{n^{4(\alpha-1)}}\int_0^{t}\int_0^\infty \left[ r\wedge(t+h-s)-r\wedge(t-s)\right]^4n^{\alpha+1}(n\theta r+1)^{-\beta}dr\,ds\nonumber\\
&=&\frac{1}{n^{4(\alpha-1)}}\int_0^{t}\int_{t-s}^{t+h-s} \left[ r-(t-s)\right]^4n^{\alpha+1}(n\theta r+1)^{-\beta}dr\,ds\nonumber\\
&&+\frac{1}{n^{4(\alpha-1)}}\int_0^{t}\int_{t+h-s}^{\infty}h^4n^{\alpha+1}(n\theta r+1)^{-\beta}dr\,ds\nonumber\\
&=&\frac{1}{n^{4(\alpha-1)}}\int_0^{t}\int_{s}^{h+s} \left( r-s\right)^4n^{\alpha+1}(n\theta r+1)^{-\beta}dr\,ds\nonumber\\
&&+\frac{1}{n^{4(\alpha-1)}}\int_0^{t}\int_{t+h-s}^{\infty}h^4n^{\alpha+1}(n\theta r+1)^{-\beta}dr\,ds.
\end{eqnarray}
 For the first term in (\ref{s-4-5}), we change the order of integration and obtain
\begin{eqnarray}\label{s-4-6}
&&\frac{1}{n^{4(\alpha-1)}}\int_0^{t}\int_{s}^{h+s} \left( r-s\right)^4n^{\alpha+1}(n\theta r+1)^{-\beta}dr\,ds\nonumber\\
&\leq&\frac{1}{n^{4(\alpha-1)}}\int_0^{h}\int_{0}^{r} \left( r-s\right)^4n^{\alpha+1}(n\theta r+1)^{-\beta}ds\,dr\nonumber\\
&&+\frac{1}{n^{4(\alpha-1)}}\int_h^{t+h}\int_{r-h}^{r} \left( r-s\right)^4n^{\alpha+1}(n\theta r+1)^{-\beta}ds\,dr\nonumber\\
&=&\frac{1}{5n^{4(\alpha-1)}}\int_0^{h}r^{5}n^{\alpha+1}(n\theta r+1)^{-\beta}dr+\frac{1}{5n^{4(\alpha-1)}}\int_h^{t+h}h^{5}n^{\alpha+1}(n\theta r+1)^{-\beta}dr\nonumber\\
&\leq&\frac{Cn^{\alpha-\beta+1}}{n^{4(\alpha-1)}}\int_0^{h} r^{5-\beta}dr+\frac{Cn^{\alpha-\beta+1}}{n^{4(\alpha-1)}}\int_{h}^{t+h}h^{5}\left(r+\frac{1}{n\theta}\right)^{-\beta}dr\nonumber\\
&\leq&Cn^{5-\beta-3\alpha}h^{6-\beta},
\end{eqnarray}
 since $\left(h+\frac{1}{n\theta}\right)^{1-\beta}-\left(t+h+\frac{1}{n\theta}\right)^{1-\beta}<\left(h+\frac{1}{n\theta}\right)^{1-\beta}<h^{1-\beta}$ for all $n\in\mathbb{N}$ and $\theta\in(0,\infty)$.
 
 For the second term in (\ref{s-4-5}), we have
\begin{eqnarray}\label{s-4-7}
\frac{1}{n^{4(\alpha-1)}}\int_0^{t}\int_{t+h-s}^{\infty}h^4n^{\alpha+1}(n\theta r+1)^{-\beta}dr\,ds
&\leq&\frac{Cn^{\alpha-\beta+1}}{n^{4(\alpha-1)}}\int_0^{t}\int_{t+h-s}^{\infty}h^4r^{-\beta}dr\,ds\nonumber\\
&=& Cn^{5-\beta-3\alpha}\int_0^{t}h^4(t+h-s)^{1-\beta}ds\nonumber \\
&\leq & Cn^{5-\beta-3\alpha}h^{6-\beta}.
\end{eqnarray}
Combining (\ref{s-4-3})-(\ref{s-4-7}), the bound (\ref{s-4-1}) follows.\\

Next we  show (\ref{s-4-2}). 
The left-hand side of \eqref{s-4-2} is bounded by
\begin{eqnarray}\label{s-4-8}
&&\frac{2}{n^{4(\alpha-1)}}\left(\int_{t}^{t+h}\int_0^\infty \left[ r\wedge(t+h-s)\right]^2n^{\alpha+1}(n\theta r+1)^{-\beta}dr\,ds\right)^{2}\nonumber\\
&+&\frac{2}{n^{4(\alpha-1)}}\left(\int_{0}^{t}\int_0^\infty \left[ r\wedge(t+h-s)-r\wedge(t-s)\right]^2n^{\alpha+1}(n\theta r+1)^{-\beta}dr\,ds\right)^{2}.
\end{eqnarray}
We bound the first term in (\ref{s-4-8}) as
\begin{eqnarray}\label{s-4-9}
&&\frac{2}{n^{4(\alpha-1)}}\left(\int_{t}^{t+h}\int_0^\infty \left[ r\wedge(t+h-s)\right]^2n^{\alpha+1}(n\theta r+1)^{-\beta}dr\,ds\right)^{2}\nonumber\\
&\leq&\frac{4}{n^{4(\alpha-1)}}\left(\int_{t}^{t+h}\int_0^{t+h-s} r^2n^{\alpha+1}(n\theta r+1)^{-\beta}dr\,ds\right)^{2}\nonumber\\
&&+\frac{4}{n^{4(\alpha-1)}}\left(\int_{t}^{t+h}\int_{t+h-s}^\infty (t+h-s)^2n^{\alpha+1}(n\theta r+1)^{-\beta}dr\,ds\right)^{2}\nonumber\\
&\leq&\frac{Cn^{2(\alpha-\beta+1)}}{n^{4(\alpha-1)}}\left(\int_{t}^{t+h}\int_0^{t+h-s} r^{2-\beta}dr\,ds\right)^{2}\nonumber\\
&&+\frac{Cn^{2(\alpha-\beta+1)}}{n^{4(\alpha-1)}}\left(\int_{t}^{t+h}\int_{t+h-s}^\infty (t+h-s)^2r^{-\beta}dr\,ds\right)^{2}\nonumber\\
&\leq&Cn^{6-2\beta-2\alpha}\left(\int_t^{t+h}(t+h-s)^{3-\beta}ds\right)^2\leq Cn^{6-2\beta-2\alpha}h^{2(4-\beta)}.
\end{eqnarray}
For the second term in (\ref{s-4-8}), by a change of variables we obtain
\begin{eqnarray}\label{s-4-10}
&&\frac{2}{n^{4(\alpha-1)}}\left(\int_{0}^{t}\int_0^\infty \left[ r\wedge(t+h-s)-r\wedge(t-s)\right]^2n^{\alpha+1}(n\theta r+1)^{-\beta}dr\,ds\right)^{2}\nonumber\\
&\leq&\frac{4}{n^{4(\alpha-1)}}\left(\int_{0}^{t}\int_{t-s}^{t+h-s} \left[ r-(t-s)\right]^2n^{\alpha+1}(n\theta r+1)^{-\beta}dr\,ds\right)^{2}\nonumber\\
&&+\frac{4}{n^{4(\alpha-1)}}\left(\int_{0}^{t}\int_{t+h-s}^{\infty} h^2n^{\alpha+1}(n\theta r+1)^{-\beta}dr\,ds\right)^{2}\nonumber\\
&=&\frac{4}{n^{4(\alpha-1)}}\left(\int_{0}^{t}\int_{s}^{h+s} (r-s)^2n^{\alpha+1}(n\theta r+1)^{-\beta}dr\,ds\right)^{2}\nonumber\\
&&+\frac{4}{n^{4(\alpha-1)}}\left(\int_{0}^{t}\int_{t+h-s}^{\infty} h^2n^{\alpha+1}(n\theta r+1)^{-\beta}dr\,ds\right)^{2}.
\end{eqnarray}
By changing the order of integration, the first term in \eqref{s-4-10} is bounded as
\begin{eqnarray}\label{s-4-11}
&&\frac{4}{n^{4(\alpha-1)}}\left(\int_{0}^{t}\int_{s}^{h+s} (r-s)^2n^{\alpha+1}(n\theta r+1)^{-\beta}dr\,ds\right)^{2}\nonumber\\
&\leq&\frac{8}{n^{4(\alpha-1)}}\left(\int_{0}^{h}\int_{0}^{r} (r-s)^2n^{\alpha+1}(n\theta r+1)^{-\beta}ds\,dr\right)^{2}\nonumber\\
&&+\frac{8}{n^{4(\alpha-1)}}\left(\int_{h}^{t+h}\int_{r-h}^{r} (r-s)^2n^{\alpha+1}(n\theta r+1)^{-\beta}ds\,dr\right)^{2}\nonumber\\
&\leq& \frac{Cn^{2(\alpha-\beta+1)}}{n^{4(\alpha-1)}}\left(\left(\int_{0}^{h}r^{3-\beta}dr\right)^{2}+\left(\int_{h}^{t+h}h^3\left(r+\frac{1}{n\theta}\right)^{-\beta}dr\right)^{2}\right)\nonumber\\
&\leq&Cn^{6-2\alpha-2\beta}h^{2(4-\beta)}.
\end{eqnarray}
The second term in (\ref{s-4-10}) can be bounded as
\begin{eqnarray}\label{s-4-12}
\frac{4n^{2(\alpha-\beta+1)}}{n^{4(\alpha-1)}}\left(\int_0^t\int_{t+h-s}^\infty h^2r^{-\beta}dr\,ds\right)^2&\leq& Cn^{6-2\alpha-2\beta}\left(\int_0^t h^2 (t+h-s)^{1-\beta}ds\right)^2\nonumber\\
&\leq&Cn^{6-2\alpha-2\beta}h^{2(4-\beta)}.
\end{eqnarray}
The bound in (\ref{s-4-2}) now follows from (\ref{s-4-8})-(\ref{s-4-11}). $\Box$

 \bibliographystyle{plain}
 \bibliography{Workload}

\end{document}